\documentclass[preprint,12pt]{elsarticle}

\usepackage{graphicx}
\usepackage{amssymb}
\usepackage{amsthm}
\usepackage{amsmath}
\usepackage{stmaryrd}
\usepackage{lineno}
\usepackage{bbold}
\usepackage{lipsum}
\usepackage{amsfonts}
\usepackage{epstopdf}
\usepackage{algorithmic}
\usepackage{cleveref}
\usepackage{accents}
\usepackage{booktabs}
\usepackage{subcaption}
\crefformat{equation}{(#2#1#3)}

\usepackage{mathrsfs}
\DeclareMathAlphabet{\mathpzc}{OT1}{pzc}{m}{it}
\usepackage{mathtools}

\usepackage{url}

\usepackage{tikz,tkz-euclide}
\usetikzlibrary{shapes.misc}
\usetikzlibrary{calc}
\tikzset{cross/.style={cross out, draw=black, minimum size=2*(#1-\pgflinewidth), inner sep=0pt, outer sep=0pt},
cross/.default={1mm}}

\newtheorem{proposition}{Proposition}

\newtheorem{remark}{Remark}

\begin{document}
\begin{frontmatter}
\date{\today}

\title{Encapsulated generalized summation-by-parts formulations for curvilinear and non-conforming meshes}


\author[add1]{Tomas Lundquist}
\author[add1]{Andrew Winters}
\author[add1,add2]{Jan Nordstr\"{o}m}

\address[add1]{Department of Mathematics, Applied Mathematics, Link\"oping University, SE-581 83 Link\"oping, Sweden}
\address[add2]{Department of Mathematics and Applied Mathematics, University of Johannesburg, P.O. Box 524, Auckland Park 2006, South Africa}  

\begin{abstract} 
We extend the construction of so-called encapsulated global summation-by-parts operators to the general case of a mesh which is not boundary conforming.
Owing to this development,  energy stable discretizations of nonlinear and variable coefficient initial boundary value problems can be formulated in simple and straightforward ways using high-order accurate operators of generalized summation-by-parts type.
Encapsulated features on a single computational block or element may include polynomial bases, tensor products as well as curvilinear coordinate transformations. Moreover, through the use of inner product preserving interpolation or projection,  the global summation-by-parts property in extended to arbitrary multi-block or multi-element meshes with non-conforming nodal interfaces.
\end{abstract}

\end{frontmatter}

\section{Introduction}

In a recent set of papers \cite{Lundquist18, aalund2019encapsulated,Lundquist22_encaps}, the idea of developing stable numerical schemes at a high level of abstraction by encapsulating low level implementation details into a global summation-by-parts (SBP) operator definition was introduced.
Such encapsulated features include curvilinear coordinate transformations on a single computational block or element, as well as non-conforming interface couplings in the case of local $h-$ or $p-$refinement.
The result is a global definition of SBP operators with diagonal norms, valid for arbitrary curvilinear multi-block or multi-element meshes. 
So far this development has focused on the classical SBP operator concept which requires a boundary conforming mesh, i.e, where the surface nodes form a subset of the volume nodes.

The primary goal of this work is to expand on this previous development by constructing encapsulated SBP formulations without the requirement of a boundary conforming mesh. For this, we use one-dimensional so-called generalized SBP operators as a baseline \cite{Fernandez14}. These operators do not require that the surface nodes are included in the volume mesh, but retain the ability to mimic the continuous energy analysis on the discrete level.  Once a set of global SBP operators on the full mesh has been constructed, numerical schemes for initial boundary value problems can be formulated in a straightforward way at a high level of abstraction, and without regards to low level implementation details such as curvilinear transformations and non-conforming interface couplings.
For nonlinear and variable coefficient problems in particular, stability proofs with the energy method can be facilitated with the introduction of additional boundary corrections in the form of penalty terms \cite{Ranocha18_gsbp}.

Extending the idea of encapsulated SBP operators for arbitrary meshes to the generalized SBP case holds several potential advantages. First,  by considering collocated pseudo-spectral operators based on the Gauss-Legendre quadrature nodes (which are not boundary conforming), a suboptimal accuracy result \cite{Lundquist18} associated with non-conforming block interfaces using classical SBP operators is avoided.  With all the details required for a stable implementation covered,  we also hope that this work may inspire the pursuit of new types of generalized SBP operators, e.g, finite difference operators, with similar favourable accuracy properties.
Finally, we anticipate that the theoretical contributions made in this paper could act as  a stepping stone to further innovations 
in key areas of numerical analysis such as embedded boundary or cut-cell methods, with the goal of simplifying mesh generation around complex geometries without sacrificing accuracy or stability.

\section{Generalized summation-by-parts operators}
In this section we define the concept of generalized SBP operators on arbitrary meshes,  and discuss the appropriate implementation leading to stable schemes for nonlinear and variable coefficient problems.

\subsection{Definition}
Consider a domain in $r-$dimensional space given by $x=(x_1,\ldots, x_r)\in\Omega\subseteq \mathbb{R}^r$, with a boundary denoted by $\partial \Omega$. Furthermore, let $\{\boldsymbol{X}_k\}_{k=1}^r$ define an arbitrary set of $N_v$ volume nodes discretizing $\Omega$, and $\{\boldsymbol{x}_k\}_{k=1}^r$ a set of $N_s$ surface nodes discretizing $\partial \Omega$. A set of high-order diagonal norm SBP partial derivative operators can in general be defined as
\begin{equation}\label{eq:SBP_def}
\mathcal{D}_{x_k} = \mathcal{P}^{-1}\mathcal{Q}_{x_k}, \quad \mathcal{Q}_{x_k}+\mathcal{Q}_{x_k}^\top  = E^\top PN_{x_k}E, \quad k = 1,\ldots , r,
\end{equation}
where $\mathcal{P}>0$ and $P>0$ are diagonal matrices containing the volume and surface quadrature weights, thus defining a discrete $L_2$ inner product on $\Omega$ and $\partial{\Omega}$, respectively. Moreover, the diagonal matrix $N_{x_k}$ contains exact or approximated values of the $k$:th outward pointing unit normal component. Finally, the rectangular matrix $E$ is a restriction or projection operator from the volume mesh to the surface mesh.

In classical definitions of SBP operators, the $N_s$ surface nodes form a subset of the $N_v$ volume nodes, meaning that the boundary restriction operator $E$ in (\ref{eq:SBP_def})
is exact; it simply consists of some permutation of $N_s$ rows from the unit matrix of dimension $N_v$.
Here we consider a generalized definition of SBP operators that do not include this requirement \cite{Fernandez14}.
In this case,  $E$ in (\ref{eq:SBP_def}) consists of a set of projection operators from the volume nodes to the surface nodes. Thus, for any smooth function $u=u(x_1,\ldots, x_r)$ defined on $\Omega$, in we consider the approximate relation
\begin{equation}\label{eq:restriction}
E  \boldsymbol{U} \approx \boldsymbol{u},
\end{equation}
where $ \boldsymbol{U}=u(\boldsymbol{X}_1, \ldots, \boldsymbol{X}_r)$ and $ \boldsymbol{u}=u (\boldsymbol{x}_1, \ldots, \boldsymbol{x}_r)$ are vectors containing the exact function values in the volume and on the surface, respectively.

One important consequence of employing inexact projections (\ref{eq:restriction}) in the SBP operator definition (\ref{eq:SBP_def}) is that the restriction operator $E$ does not commute with respect to elementwise multiplication. Consider two smooth mesh functions $\boldsymbol{\Phi}$ and $\boldsymbol{U}$. For the product between these two we can write, due to (\ref{eq:restriction}),
\begin{equation}\label{eq:restriction_inexact}
E(\underline{\boldsymbol{\Phi}}\boldsymbol{U} ) \approx \underline{ \boldsymbol{\phi}}\boldsymbol{u},  \quad  \boldsymbol{\phi} = E\boldsymbol{\Phi}, \quad \boldsymbol{u} =E\boldsymbol{U},
\end{equation}
where $\underline{\boldsymbol{\Phi}}$ and $\underline{\boldsymbol{\phi}}$  denote the diagonal matrices obtained by inserting the coefficient vectors $\boldsymbol{\Phi}$ and $\boldsymbol{\phi}$ along the diagonal.  Due to the lack of strict equality in (\ref{eq:restriction_inexact}) for arbitrary vectors $\boldsymbol{\Phi}$ and $\boldsymbol{U}$, it was shown in \cite{Nordstrom17_gsbp} that a naive implementation of generalized SBP operators in the context of nonlinear and variable coefficient problems is not conducive to stability proofs with the energy method. However, it was later demonstrated in \cite{Ranocha18_gsbp} that the issue identified in \cite{Nordstrom17_gsbp} can be resolved by the introduction of boundary corrections in the form of penalty terms.

 \subsection{The energy method}
The fundamental building block for analyzing continuous initial boundary value problems with the energy method is integration-by-parts.
In the context of general nonlinear and variable coefficient problems,  we consider multi-dimensional integration-by-parts rules of the type
\begin{equation}\label{eq:IBP_varcoef}
\big(v, \psi(\phi u)_{x_k}\big)_\Omega + \big(u, \phi (\psi v)_{x_k} \big)_\Omega =  \oint_{\partial \Omega}(\phi u\psi v) n_{x_k}, \quad k=1,\ldots, r.
\end{equation}
where $\phi$, $\psi$, $u$, $v$ are smooth scalar functions over $\Omega$, and $n_{x_k}$ denotes the $k$:th continuous component of the outward pointing unit normal vector.
A straightforward (naive) implementation of the $x_k-$derivative SBP operator $\mathcal{D}_{x_k}= \mathcal{P}^{-1}\mathcal{Q}_{x_k}$ in (\ref{eq:SBP_def}) leads to the discrete approximation
\begin{equation}\label{eq:varcoef_naive}
\boldsymbol{\Psi}(\boldsymbol{\Phi} \boldsymbol{U})_{x_k} \approx  
 \underline{\boldsymbol{\Psi}} \mathcal{D}_{x_k}
(\underline{\boldsymbol{\Phi}} \boldsymbol{U})= 
 \mathcal{P}^{-1} (\underline{\boldsymbol{\Psi}}\mathcal{Q}_{x_k}\underline{\boldsymbol{\Phi}}) \boldsymbol{U},
\end{equation}
where in the last step we have used the fact that $\mathcal{P}^{-1}$ and $\boldsymbol{\Psi}$ are diagonal matrices, and thus commute.
From the SBP property (\ref{eq:SBP_def}), we have
\begin{equation}\label{eq:SBP_varcoef_naive}
\underline{\boldsymbol{\Psi}}\mathcal{Q}_{x_k}\underline{\boldsymbol{\Phi}}+(\underline{\boldsymbol{\Phi}}\mathcal{Q}_{x_k}\underline{\boldsymbol{\Psi}})^\top  =   (E\underline{\boldsymbol{\Psi}})^\top  PN_{ x_k}(E\underline{\boldsymbol{\Phi}}).
\end{equation}
We thus get, using (\ref{eq:restriction_inexact}),
\begin{equation}\label{eq:SBP_vec_naive}
\big(\boldsymbol{V},  \underline{\boldsymbol{\Psi}} \mathcal{D}_{x_k}
(\underline{\boldsymbol{\Phi}} \boldsymbol{U}) \big)_{\mathcal{P}}+ \big(\boldsymbol{U},   \underline{\boldsymbol{\Phi}} \mathcal{D}_{x_k}
(\underline{\boldsymbol{\Psi}} \boldsymbol{V}) \big)_{\mathcal{P}}
=(E\underline{\boldsymbol{\Psi}}\boldsymbol{V})^\top  PN_{ x_k}(E\underline{\boldsymbol{\Phi}}\boldsymbol{U})
\approx  (\underline{\boldsymbol{\psi}}\boldsymbol{v})^\top PN_{x_k}(\underline{\boldsymbol{\phi}}\boldsymbol{u})  ,
\end{equation}
which \textit{approximately} mimics integration-by-parts (\ref{eq:IBP_varcoef}) discretely in an inner product sense.  In the special case of classical SBP operators, 
this inner product relation is satisfied with \textit{equality} using arbitrary vectors $\boldsymbol{\Phi}$, $\boldsymbol{\Psi}$, $\boldsymbol{U}$, $\boldsymbol{V}$ of length $N_v$, which enables stability proofs with the energy method. In order to obtain the same result in the generalized SBP case,  a modification to the straightforward approximation in (\ref{eq:varcoef_naive}) is required.  Following a similar approach as in \cite{Ranocha18_gsbp}, we will do this by applying a  boundary correction penalty term. 

We start by adding and subtracting the desired exact boundary term from (\ref{eq:SBP_varcoef_naive}),
\begin{equation}\label{eq:step1}
\underline{\boldsymbol{\Psi}}\mathcal{Q}_{x_k}\underline{\boldsymbol{\Phi}}+(\underline{\boldsymbol{\Phi}}\mathcal{Q}_{x_k}\underline{\boldsymbol{\Psi}})^\top  = (\underline{ \boldsymbol{\psi}}E)^\top PN_{x_k}(\underline{\boldsymbol{\phi}}E)-(\underline{ \boldsymbol{\psi}}E)^\top PN_{x_k}(\underline{\boldsymbol{\phi}}E)+  (E\underline{\boldsymbol{\Psi}})^\top  PN_{ x_k}(E\underline{\boldsymbol{\Phi}}) .
\end{equation}
The remaining (undesired) part of the right hand side above can be expanded according to
\begin{equation}\label{eq:step2}
\begin{aligned}
& -(\underline{ \boldsymbol{\psi}}E)^\top PN_{x_k}(\underline{\boldsymbol{\phi}}E)+(E\underline{\boldsymbol{\Psi}})^\top  PN_{ x_k}(E\underline{\boldsymbol{\Phi}})  \\
= \ & -\frac{1}{2}\big(\underline{\boldsymbol{\psi}}E+E\underline{\boldsymbol{\Psi}} \big)^\top PN_{ x_k}(\underline{\boldsymbol{\phi}}E-E\underline{\boldsymbol{\Phi}}) - \Big[\frac{1}{2}\big(\underline{\boldsymbol{\phi}}E+E\underline{\boldsymbol{\Phi}} \big)^\top PN_{ x_k}(\underline{\boldsymbol{\psi}}E-E\underline{\boldsymbol{\Psi}})\Big]^\top .
\end{aligned}
\end{equation}
From this, we identify a suitable penalty correction to the naive approximation in (\ref{eq:varcoef_naive}) as follows,
\begin{equation}\label{eq:varcoef_disc}
\boldsymbol{\Psi}(\boldsymbol{\Phi} \boldsymbol{U})_{x_k} \approx  \mathcal{D}_{x_k,\boldsymbol{\Phi}, \boldsymbol{\Psi}}\boldsymbol{U} = \mathcal{P}^{-1}\mathcal{Q}_{x_k,\boldsymbol{\Phi}, \boldsymbol{\Psi}}\boldsymbol{U},
\end{equation}
where the modified undivided difference matrix $\mathcal{Q}_{x_k,\boldsymbol{\Phi}, \boldsymbol{\Psi}}$ contains the volume contribution from (\ref{eq:varcoef_naive}) as well as an added boundary correction term,
\begin{equation}\label{eq:varcoef_Q}
\mathcal{Q}_{x_k,\boldsymbol{\Phi}, \boldsymbol{\Psi}}=
\underline{\boldsymbol{
\Psi}}\mathcal{Q}_{x_k}\underline{\boldsymbol{\Phi}} +\frac{1}{2}\big(\underline{\boldsymbol{\psi}}E+E\underline{\boldsymbol{\Psi}} \big)^\top 
PN_{ x_k}(\underline{\boldsymbol{\phi}}E-E\underline{\boldsymbol{\Phi}}).
\end{equation}
We have thus proven the first main result of this paper.
\begin{proposition}\label{thm:correction}
Consider the discrete variable coefficient approximation (\ref{eq:varcoef_disc}), (\ref{eq:varcoef_Q}), where $\mathcal{Q}_{x_k}$ satisfies the partial derivative SBP property (\ref{eq:SBP_def}).
The following variable coefficient SBP property is then satisfied for arbitrary vectors $\boldsymbol{\Phi}$ and $\boldsymbol{\Psi}$ of length $N_v$,
\begin{equation}\label{eq:SBP_varcoef}
\mathcal{Q}_{x_k,\boldsymbol{\Phi}, \boldsymbol{\Psi}}+\mathcal{Q}_{x_k,\boldsymbol{\Psi}, \boldsymbol{\Phi}}^T=(\underline{ \boldsymbol{\psi}}E)^\top PN_{x_k}(\underline{\boldsymbol{\phi}}E).
\end{equation}
\begin{proof}
The result follows directly from expanding $\mathcal{Q}_{x_k,\boldsymbol{\Phi}, \boldsymbol{\Psi}}+\mathcal{Q}_{x_k,\boldsymbol{\Psi}, \boldsymbol{\Phi}}^T$ using the definition in (\ref{eq:varcoef_Q}), and then inserting (\ref{eq:step1}) followed by (\ref{eq:step2}).
\end{proof}
\end{proposition}
With (\ref{eq:SBP_varcoef}), in place,  the inner product statement (\ref{eq:SBP_vec_naive}) is modified to become
\begin{equation}\label{eq:SBP_vec}
\big(\boldsymbol{V},  \mathcal{D}_{x_k,\boldsymbol{\Phi}, \boldsymbol{\Psi}}\boldsymbol{U} \big)_{\mathcal{P}}+ \big(\boldsymbol{U},   \mathcal{D}_{x_k,\boldsymbol{\Psi}, \boldsymbol{\Phi}} \boldsymbol{V} \big)_{\mathcal{P}}
=(\underline{ \boldsymbol{\psi}}E\boldsymbol{V})^\top  PN_{ x_k}(\underline{ \boldsymbol{\phi}}E\boldsymbol{U})
=  (\underline{\boldsymbol{\psi}}\boldsymbol{v})^\top PN_{x_k}(\underline{\boldsymbol{\phi}}\boldsymbol{u})  ,
\end{equation}
which exactly mimics integration-by-parts (\ref{eq:IBP_varcoef}) discretely in an inner product sense.
\begin{remark}
The discrete boundary quantities $\boldsymbol{\phi}$, $\boldsymbol{\psi}$ in (\ref{eq:varcoef_Q}) are not formally required to satisfy the exact relations $\boldsymbol{\phi} = E\boldsymbol{\Phi}$ and $\boldsymbol{\psi} =E\boldsymbol{\Psi}$. In a nonlinear setting however, this is a convenient choice since the exact values of the continuous functions $\phi$ and $\psi$ may be solution dependent, and thus unavailable.
\end{remark}
With the discrete energy method in place, what remains is to demonstrate the validity and versatility of the encapsulated partial derivative SBP property (\ref{eq:SBP_def}) on complex, multi-domain geometries.

\section{Polynomial bases, tensor products and curvilinear transformations}
\label{sec:polytens}
It is well known that the use of SBP operators defined in one space dimension can be extended to general coordinates in multiple dimensions with a combination of tensor (or Kronecker) products as well as curvilinear coordinate transformations. In this section we demonstrate how these techniques can be incorporated into the encapsulated SBP operator definition (\ref{eq:SBP_def}), thus extending the methodology previously outlined in \cite{Lundquist18,aalund2019encapsulated} from classical to generalized SBP operators.

We first note that in one spatial dimension, the SBP definition (\ref{eq:SBP_def}) on an interval $x\in\Omega=[a,b]$ simplifies to
\begin{equation}\label{eq:SBP_1D}
D_x = P_x^{-1}Q_x, \quad Q_x+Q_x^\top  = 
\begin{pmatrix}
e_a & e_b
\end{pmatrix}
\begin{pmatrix}
-1 \\ & 1
\end{pmatrix}
\begin{pmatrix}
e_a^\top  \\ e_b^\top 
\end{pmatrix}=
e_b e_b^\top -e_ae_a^\top  ,
\end{equation}
where $e_a=e_{x=a}$ and $e_b=e_{x=b}$ are column vectors such that $e_a^\top $, $e_b^\top $ defines a pair of projection operators to the left ($x=a$) and the right ($x=b$) boundary node, respectively.

\subsection{Polynomial bases}\label{sec:poly}
Although other types of numerical methods (such as finite difference methods) could be considered, the main target with the present study is SBP operators based on pseudo-spectral collocation.  The fundamental link between finite element approximations and dense norm SBP operators was discussed already in the seminal work \cite{ref:KREI74} as well as later, e.g, in \cite{carp96,nordstrom2006conservative}. More importantly for our purposes however, by collocating a spectral (polynomial) finite element first derivative approximation at the quadrature nodes,  a diagonal norm SBP operator is obtained \cite{gassner13, Fernandez14}.

On a single pseudo-spectral element in $1$D, the diagonal norm matrix $P_x$ in (\ref{eq:SBP_1D}) is given by the collocation quadrature weights, while $Q_x$ is given by the exact integral
\begin{equation}\label{eq:pseudo-op}
Q_x = \int_\Omega \mathcal{L}\mathcal{L}_x^\top  dx,
\end{equation}
where the vector function $\mathcal{L}=\mathcal{L}(x)$ contains the set of interpolating Lagrange polynomials. Using integration-by-parts, the formulation in (\ref{eq:pseudo-op}) leads to
\begin{equation*}
Q_x + Q_x^\top  = \mathcal{L}(b)\mathcal{L}(b)^\top -\mathcal{L}(a)\mathcal{L}(a)^\top  =  e_b e_b^\top -e_ae_a^\top ,
\end{equation*}
i.e, (\ref{eq:SBP_1D}) is indeed satisfied.  
Note that if $a$ and $b$ are included in the set of collocation quadrature nodes, such as for Gauss-Lobatto quadrature rules,  then the Kronecker delta property of Lagrange interpolating polynomials implies that $e_a^\top =\mathcal{L}(a)^\top $ and $e_b^\top =\mathcal{L}(b)^\top $, i.e, the boundary projection operators reduce to the first and last unit matrix rows. Thus,  the operator $D_x$ is SBP in the classical sense \cite{gassner13}.  If the quadrature nodes are not boundary conforming, as in the case of the optimally accurate Gauss-Legendre rules, then the accuracy of $e_a^\top $ and $e_b^\top $ is determined by the polynomial degree $N$.

\subsubsection{Tensor products}
For example, starting from the $1$D operators $D_\xi=P_ \xi^{-1}Q_\xi$ and $D_\eta=P_\eta^{-1}Q_\eta$ satisfying (\ref{eq:SBP_1D}) on the unit interval $[0,1]$,
we define a corresponding pair of $2$D operators on the reference unit square $(\xi,\eta)\in\hat{\Omega}=[0,1]^2$, given by $\mathcal{D}_\xi=D_\xi\otimes I_\eta$ and $\mathcal{D}_\eta=I_\xi\otimes D_\eta $. We express these as $\mathcal{D}_\xi=\hat{\mathcal{P}}^{-1}\mathcal{Q}_\xi$ and $\mathcal{D}_\eta= \hat{\mathcal{P}}^{-1}\mathcal{Q}_\eta$, where
\begin{equation}\label{eq:SBP_tensor}
\hat{\mathcal{P}} = (P_\xi\otimes P_\eta), \quad \mathcal{Q}_\xi = (Q_\xi \otimes 
P_\eta), \quad \mathcal{Q}_\eta = (P_\xi\otimes Q_\eta).
\end{equation}
It is straightforward to verify that the SBP property (\ref{eq:SBP_def}) is satisfied with respect to the block matrices
\begin{equation}\label{eq:tensor_boun_1}
 E =  
\begin{pmatrix}
e_{\xi=0}^\top \otimes I_\eta  \\[0.1cm]
e_{\xi=1}^\top \otimes I_y \\[0.1cm]
I_\eta\otimes e_{\eta=0}^\top  \\[0.1cm]
I_\eta\otimes e_{\eta=1}^\top 
\end{pmatrix},\quad \hat{P} = \begin{pmatrix}
P_\eta \\ & P_\eta \\ && P_\xi \\ &&& P_\xi
\end{pmatrix}, 
\end{equation}
as well as
\begin{equation}\label{eq:tensor_boun_2}
 N_\xi =  \begin{pmatrix}
I_\eta  \\ & -I_\eta \\ && 0 \\ &&& 0
\end{pmatrix},  \quad  N_\eta  = \begin{pmatrix}
0  \\ & 0 \\ && I_\xi \\ &&& -I_\xi
\end{pmatrix},
\end{equation}
with the blocks representing the East, West, North and South faces to $\hat{\Omega}$. 
There is no need here to distinguish between classical and generalized SBP operators, the tensor product formulation is identical for either type of operator.

\subsubsection{Curvilinear coordinates}
\label{sec:curvi}
For simplicity of presentation we focus on the $2$D case, although the analysis can be extended to an arbitrary number of space dimensions. Thus, consider a $2$D domain given by $(x,y)\in \Omega$, where the functions $x=x(\xi,\eta)$ and $y=y(\xi, \eta)$ define a curvilinear transformation from the reference unit square $(\xi,\eta)\in \hat{\Omega}=[0,1]^2$. 
Applying the chain rule, the partial $x-$ and $y-$derivatives of a function $\phi$ on can be expressed on the skew-symmetric form
\begin{equation}\label{eq:curvi_cont}
\begin{aligned}
\phi_x=\ & \frac{1}{2} J^{-1}\Big[ (y_\eta \phi)_\xi+y_\eta \phi_\xi-(y_\xi \phi)_\eta-y_\xi \phi_\eta\Big] \\
\phi_y=\ & \frac{1}{2} J^{-1}\Big[ (x_\xi \phi)_\eta+x_\xi \phi_\eta-(x_\eta \phi)_\xi-x_\eta \phi_\xi\Big],
\end{aligned}
\end{equation}
where $J=x_\xi y_\eta - y_\xi x_\eta $ denotes the Jacobian determinant of the transformation.  

On the discrete side, we start from a pair of SBP operators $\mathcal{D}_\xi=\hat{\mathcal{P}}^{-1}\mathcal{Q}_\xi$ and $\mathcal{D}_\eta=\hat{\mathcal{P}}^{-1}\mathcal{Q}_\eta$ on the reference domain $\hat{\Omega}$. A straightforward (naive) approximation of (\ref{eq:curvi_cont}) is given by $\boldsymbol{\Phi}_x\approx \mathcal{P}^{-1}\tilde{\mathcal{Q}}_{x}\boldsymbol{\Phi}$ and $\boldsymbol{\Phi}_y\approx \mathcal{P}^{-1}\tilde{\mathcal{Q}}_{y}\boldsymbol{\Phi}$, where $\mathcal{P} = \underline{\boldsymbol{J}}\hat{\mathcal{P}}$ contains the rescaled quadrature weights, and
\begin{equation}\label{eq:curvi_Q_tilde}
\begin{aligned}
\tilde{\mathcal{Q}}_x = \ & \frac{1}{2} \big(\mathcal{Q}_\xi\underline{\boldsymbol{Y}_\eta} +\underline{\boldsymbol{Y}_\eta}\mathcal{Q}_\xi -\mathcal{Q}_{\eta}\underline{\boldsymbol{Y}_\xi}-\underline{\boldsymbol{Y}_\xi}\mathcal{Q}_{\eta}\big) \\
 \tilde{\mathcal{Q}}_y = \ & \frac{1}{2}\big(\mathcal{Q}_{\eta}\underline{\boldsymbol{X}_\xi}+\underline{\boldsymbol{X}_\xi}\mathcal{D}_{\eta}- \mathcal{Q}_\xi\underline{\boldsymbol{X}_\eta} - \underline{\boldsymbol{X}_\eta}\mathcal{Q}_\xi \big).
 \end{aligned}
\end{equation}
Indeed, in the classical SBP operator case, this formulation corresponds to the encapsulated operator definition considered in \cite{Lundquist18,aalund2019encapsulated}.

Next, in the generalized SBP case, we can apply Proposition \ref{thm:correction} on the level of the reference domain $ \hat{\Omega}$ in order to derive an appropriate penalty correction to (\ref{eq:curvi_Q_tilde}).  In particular, consider the definition
\begin{equation}\label{eq:curvi_disc_corr}
\begin{aligned}
\mathcal{Q}_x = \ & \frac{1}{2} \big(\mathcal{Q}_{\xi, \boldsymbol{Y}_\eta, \mathbb{1}}+\mathcal{Q}_{\xi, \mathbb{1}, \boldsymbol{Y}_\eta} -\mathcal{Q}_{\eta, \boldsymbol{Y}_\xi, \mathbb{1}}-\mathcal{Q}_{\eta, \mathbb{1},\boldsymbol{Y}_\xi}\big) \\
\mathcal{Q}_y = \ & \frac{1}{2} \big(\mathcal{Q}_{\eta, \boldsymbol{X}_\xi, \mathbb{1}}+\mathcal{Q}_{\eta, \mathbb{1}, \boldsymbol{X}_\xi} -\mathcal{Q}_{\xi, \boldsymbol{X}_\eta, \mathbb{1}}-\mathcal{Q}_{\xi, \mathbb{1},\boldsymbol{X}_\eta}\big) ,
 \end{aligned}
\end{equation}
where each variable coefficient term is defined as in (\ref{eq:varcoef_Q}).  Note that for the vectors $\mathbb{1}$ and $\boldsymbol{1}$ containing the constant value $1$,
we have $ \underline{\boldsymbol{1}}E+E\underline{\mathbb{1}}=2E$ as well as $\underline{\boldsymbol{1}}E-E\underline{\mathbb{1}}=0$.
Therefore,  starting from the general definition in (\ref{eq:varcoef_Q}), the terms in (\ref{eq:curvi_Q_tilde}) are simplified to yield
\begin{equation}\label{eq:curvi_varcoef}
\begin{aligned}
\mathcal{Q}_{\xi, \boldsymbol{Y}_\eta, \mathbb{1}} = \ & \mathcal{Q}_\xi \underline{\boldsymbol{Y}_\eta} + E^\top \hat{P}N_{\xi} \big(\underline{\boldsymbol{y}_\eta}E-E\underline{\boldsymbol{Y}_\eta}\big)\quad \mathcal{Q}_{\xi, \mathbb{1},\boldsymbol{Y}_\eta} =  \underline{\boldsymbol{Y}_\eta} \mathcal{Q}_\xi \\
\mathcal{Q}_{\eta, \boldsymbol{Y}_\xi, \mathbb{1}} = \ & \mathcal{Q}_\eta \underline{\boldsymbol{Y}_\xi} + E^\top \hat{P}N_{\eta} \big(\underline{\boldsymbol{y}_\xi}E-E\underline{\boldsymbol{Y}_\xi}\big)\quad \mathcal{Q}_{\eta, \mathbb{1},\boldsymbol{Y}_\xi} =  \underline{\boldsymbol{Y}_\xi} \mathcal{Q}_\eta \\
\mathcal{Q}_{\eta, \boldsymbol{X}_\xi, \mathbb{1}} = \ & \mathcal{Q}_\eta \underline{\boldsymbol{X}_\xi} + E^\top \hat{P}N_{\eta} \big(\underline{\boldsymbol{x}_\xi}E-E\underline{\boldsymbol{x}_\xi}\big)\quad \mathcal{Q}_{\eta, \mathbb{1},\boldsymbol{X}_\xi} =  \underline{\boldsymbol{X}_\xi} \mathcal{Q}_\eta \\
\mathcal{Q}_{\xi, \boldsymbol{X}_\eta, \mathbb{1}} = \ & \mathcal{Q}_\xi \underline{\boldsymbol{X}_\eta} + E^\top \hat{P}N_{\xi} \big(\underline{\boldsymbol{x}_\eta}E-E\underline{\boldsymbol{x}_\eta}\big)\quad \mathcal{Q}_{\xi, \mathbb{1},\boldsymbol{X}_\eta} =  \underline{\boldsymbol{X}_\eta} \mathcal{Q}_\xi .
 \end{aligned}
\end{equation}
where discrete volume and surface values of the metric coefficients $x_\xi$, $x_\eta$, $y_\xi$, $y_\eta$ can be either exact or approximated.
Written in terms of a penalty correction to the naive formulation in (\ref{eq:curvi_Q_tilde}), the operators $\mathcal{Q}_x$ and $\mathcal{Q}_y$ in (\ref{eq:curvi_disc_corr}) thus become
\begin{equation}\label{eq:curvi_disc_x}
\begin{aligned}
\mathcal{Q}_x = \ & \tilde{\mathcal{Q}}_x+\frac{1}{2}E^\top \hat{P}\Big[N_\xi\big(\underline{\boldsymbol{y}_\eta}E-E\underline{\boldsymbol{Y}_\eta}\big)-N_{\eta} \big(\underline{\boldsymbol{y}_\xi}E-E\underline{\boldsymbol{Y}_\xi}\big)\Big] \\
\mathcal{Q}_y = \ & \tilde{\mathcal{Q}}_y+\frac{1}{2}E^\top \hat{P}\Big[N_\eta\big(\underline{\boldsymbol{x}_\xi}E-E\underline{\boldsymbol{X}_\xi}\big)-N_{\xi} \big(\underline{\boldsymbol{x}_\eta}E-E\underline{\boldsymbol{X}_\eta}\big)\Big] .
\end{aligned}
\end{equation}

Next, we verify that $\mathcal{Q}_x$ and $\mathcal{Q}_y$ in (\ref{eq:curvi_disc_x})  satisfy the encapsulated SBP property (\ref{eq:SBP_def}) on the level of the physical domain $\Omega$. We first note as in \cite{Lundquist18} that the stretching factor between $\partial\hat{\Omega}$ and $\partial \Omega$ on the four faces is given by the following continuous function,
\begin{equation*}
j (x,y)=     \begin{cases}
      \sqrt{x_\eta^2+y_\eta^2}, & \text{East, West.}\\
      \sqrt{x_\xi^2+y_\xi^2}, & \text{North, South}.
    \end{cases}  
\end{equation*}
It is therefore natural to define a rescaled quadrature rule on $\partial \Omega$ as
\begin{equation}\label{eq:curvi_boun_P} 
P=\underline{\boldsymbol{j}}\hat{P},
\end{equation}
where $\boldsymbol{j}$ contains exact or approximated nodal values of this function $j$.
Moreover, the outward pointing unit normal $n=(n_x, \ n_y)^\top $ on $\partial \Omega$ is given by
\begin{equation*}
jn_x = \begin{cases}
\phantom{-}y_\eta , & \text{East}\\
      -y_\eta , & \text{West}\\
-y_\xi , & \text{North}\\
\phantom{-}y_\xi, & \text{South.}
    \end{cases}  , \quad
    jn_y = \begin{cases}
      -x_\eta, & \text{East}\\
      \phantom{-}x_\eta, & \text{West}\\
\phantom{-}x_\xi, & \text{North}\\
      -x_\xi, & \text{South.}
    \end{cases}  
\end{equation*}
Using (\ref{eq:tensor_boun_2}), we thus define
\begin{equation}\label{eq:curvi_boun}
N_x = \underline{\boldsymbol{j}}^{-1}(N_\xi\underline{\boldsymbol{y}_\eta} - N_\eta\underline{\boldsymbol{y}_\xi}), \quad N_y = \underline{\boldsymbol{j}}^{-1}(N_\eta\underline{\boldsymbol{x}_\xi} - N_\xi\underline{\boldsymbol{x}_\eta}).
\end{equation}
We are now ready to prove the second main result of this paper.
\begin{proposition}\label{thm:curvi_SBP}
The curvilinear operators defined in (\ref{eq:curvi_disc_x}) satisfy the encapsulated physical domain SBP property (\ref{eq:SBP_def}) with respect to the discrete surface quadrature and normal vector components defined in (\ref{eq:curvi_boun_P}), (\ref{eq:curvi_boun}).

\begin{proof}
We prove the result for the $x-$coordinate case as the $y-$coordinate case is done analogously. By applying Proposition \ref{thm:correction} to the first two lines in (\ref{eq:curvi_varcoef}), it directly follows that $\mathcal{Q}_x$ satisfies
\begin{equation*}
\mathcal{Q}_x+\mathcal{Q}_x^\top  = E^\top \hat{P} ( N_\xi\underline{\boldsymbol{y}_\eta} -N_\eta\underline{\boldsymbol{y}_\xi})E .
\end{equation*}
By inserting $I=\underline{\boldsymbol{j}} \ \underline{\boldsymbol{j}}^{-1}$ above,  the SBP property (\ref{eq:SBP_def}) now follows from the definitions in (\ref{eq:curvi_boun_P}), (\ref{eq:curvi_boun}).
\end{proof}
\end{proposition}
We further prove the following result regarding the consistency of the encapsulated SBP operator formulation (\ref{eq:curvi_disc_x}).
\begin{proposition}\label{thm:consistency}
Assume that the following consistency conditions are satisfied by the reference domain operators,
\begin{equation}\label{eq:consistency_ref}
\mathcal{Q}_\xi\mathbb{1} = 0, \quad \mathcal{Q}_\eta \mathbb{1} = 0,
\end{equation}
as well as, for the boundary projection operator,
\begin{equation} \label{eq:consistency_project}
E\mathbb{1}=\boldsymbol{1},
\end{equation}
where $\mathbb{1}$ and $\boldsymbol{1}$ are vectors with the constant value $1$ everywhere.
Then the physical domain operators in (\ref{eq:curvi_disc_x}) are also consistent, provided that the following three additional assumptions hold,
\begin{enumerate}[(i)]
\item The metric coefficients on the volume mesh are approximated according to $\boldsymbol{X}_\xi=\mathcal{D}_\xi\boldsymbol{X}$, $\boldsymbol{X}_\eta=\mathcal{D}_\eta\boldsymbol{X}$, $\boldsymbol{Y}_\xi=\mathcal{D}_\xi\boldsymbol{Y}$, $\boldsymbol{Y}_\eta=\mathcal{D}_\eta\boldsymbol{Y}$.
\item The metric coefficients on the surface mesh are given by the projections $\boldsymbol{x}_\xi = E\boldsymbol{X}_\xi$, $\boldsymbol{x}_\eta = E\boldsymbol{X}_\eta$, $\boldsymbol{y}_\xi = E\boldsymbol{Y}_\xi$, $\boldsymbol{y}_\eta = E\boldsymbol{Y}_\eta$.
\item The two operators $\mathcal{D}_\xi$, $\mathcal{D}_\eta$ commute, i.e, $\mathcal{D}_\xi\mathcal{D}_\eta = \mathcal{D}_\eta\mathcal{D}_\xi$. 
\end{enumerate}
\begin{proof}
Again, we prove the result for the $x-$coordinate case. For an arbitrary vector $\boldsymbol{\phi}$, we can first write $\underline{\boldsymbol{\phi}}E\mathbb{1}=\boldsymbol{\phi}$ due to (\ref{eq:consistency_project}). Thefore,  multiplying $\mathcal{Q}_x$ in (\ref{eq:curvi_disc_x}) with the constant vector $\mathbb{1}$ yields
\begin{equation*}
\mathcal{Q}_x \mathbb1 =\tilde{\mathcal{Q}}_x\mathbb{1}+\frac{1}{2}E^\top \hat{P}\Big[N_\xi 
(  \boldsymbol{y}_\eta - E\boldsymbol{Y}_\eta  ) - N_\eta (  \boldsymbol{y}_\xi- E\boldsymbol{Y}_\xi  ) \Big].
\end{equation*}
Since $\boldsymbol{y}_\eta = E\boldsymbol{Y}_\eta$ and $\boldsymbol{y}_\xi = E\boldsymbol{Y}_\xi$ holds by assumption (ii), this yields $\mathcal{Q}_x \mathbb1 =\tilde{\mathcal{Q}}_x\mathbb{1}$. Next, using (\ref{eq:consistency_ref}), we can write,
\begin{equation*}
\tilde{\mathcal{Q}}_x\mathbb1  = \frac{1}{2} \big(\mathcal{Q}_\xi\boldsymbol{Y}_\eta -\mathcal{Q}_{\eta}\boldsymbol{Y}_\xi\big) = \frac{1}{2}\mathcal{P}\big(\mathcal{D}_{\xi}\mathcal{D}_{\eta}-\mathcal{D}_{\eta}\mathcal{D}_{\xi}\big)\boldsymbol{Y},
\end{equation*}
where in the last step we have used assumption (i).
Thus, since $\mathcal{D}_\xi\mathcal{D}_\eta = \mathcal{D}_\eta\mathcal{D}_\xi$ holds by assumption (iii), the consistency condition $\mathcal{Q}_x \mathbb1 =0$ follows.
\end{proof} 
\end{proposition}
\begin{remark}
If tensor products are employed as in (\ref{eq:SBP_tensor}) in order to construct the reference domain operators $\mathcal{D}_\xi$ and $\mathcal{D}_\eta$, then these operators automatically satisfy the commutative property $\mathcal{D}_\xi\mathcal{D}_\eta = \mathcal{D}_\eta\mathcal{D}_\xi$ required in Proposition \ref{thm:consistency}. Furthermore, they inherit the consistency properties (\ref{eq:consistency_ref}) and (\ref{eq:consistency_project}) from the one-dimensional operators $D_\xi$ and $D_\eta$.
\end{remark}

\section{Multi-block/element meshes}
Next, we discuss the construction and implementation of encapsulated generalized SBP operators (\ref{eq:SBP_def}) on arbitrary multi-block or multi-element meshes with non-conforming nodal interfaces.
For this purpose, the same procedure from \cite{Lundquist22_encaps} focusing on classical SBP operators applies without modification to the generalized SBP case well. 
To maintain maximal generality,
we here discuss the "many-to-many" block or element coupling case, thereby extending the more simple one-to-one coupling procedure previously considered in \cite{Lundquist22_encaps}.  In other words, we make no assumptions regarding either the shape or the topology of the multi-block/element geometry.

Thus, consider a generic partitioning $\Omega = \bigcup_{m=1}^{M}\Omega_m$ of the spatial domain into $M$ non-overlapping subdomains,  and let $\{\mathcal{D}_{x_k,m}
\}_{m=1}^{M}$ denote a set of SBP operators (\ref{eq:SBP_def}) approximating the $x_k$-derivative on each $\Omega_m$.
We split each subdomain boundary into two parts, $\partial \Omega_m=\partial^e \Omega_m \cup \partial^i\Omega_m$, where $\partial^e \Omega_m$ is located on the exterior boundary to $\Omega$, while $\partial ^i\Omega_m$ forms internal interfaces with neighbouring subdomains. Note that $\partial^e\Omega_m$ is an empty set if $\Omega_m$ is completely located in the domain interior, and surrounded by other subdomains in all directions.
In (\ref{eq:SBP_def}), we accordingly separate the right-hand-side contributions into two parts,
\begin{equation}\label{eq:SBP_split}
\begin{aligned}
\mathcal{Q}_{x_k,m}+\mathcal{Q}_{x_k,m}^\top =E_m^\top P_mN_{x_k,m}E_m 
= \ & (E_m^e)^\top P_m^e N_{x_k,m}^e E_m^e + (E_m^e)^\top P_m^iN_{x_k,m}^i E_m^i,
\end{aligned}
\end{equation}
where
\begin{equation*}
E_m = \begin{pmatrix}
E_m^e \\E_m^i
\end{pmatrix}, \quad 
P_m = \begin{pmatrix}
P_m^e \\ & P_m^i
\end{pmatrix}, \quad 
N_{x_k,m}= \begin{pmatrix}
N_{x_k,m}^e \\ & N_{x_k,m}^i
\end{pmatrix}.
\end{equation*}
Furthermore, to assemble all subdomain contributions into a set of global block matrices we, for notational convenience,  introduce three sets of additional restriction operators, 
\begin{equation}\label{eq:restr_ops}
\begin{aligned}
\tilde{I}_m^v=  \ & \begin{pmatrix}
0 & \ldots &  0 & I_m^v & 0 & \ldots & 0
\end{pmatrix} \\
\tilde{I}_m^e =  \ & \begin{pmatrix}
0 & \ldots &  0 & I_m^e & 0 & \ldots & 0
\end{pmatrix} \\ 
\tilde{I}_m^i =  \ & \begin{pmatrix}
0 & \ldots &  0 & I_m^i & 0 & \ldots & 0
\end{pmatrix}, \\ 
\end{aligned}
\end{equation}
for $m=1,\ldots, M$, where $I_m^v$, $I_m^e$ and $I_m^i $ are unit matrices corresponding to the number of volume,  boundary and interface nodes on the single subdomain $\Omega_m$, respectively.
In other words, $\tilde{I}_m^v$ in (\ref{eq:restr_ops}) acts to restrict from the full volume mesh on $\Omega$ onto $\Omega_m$, while $\tilde{I}_m^e$ and $\tilde{I}_m^i $ restrict from the full set of exterior boundary and interior interface nodes to only those on $\partial^e \Omega_m$ and $\partial^i \Omega_m$, respectively. 

Next, we assemble all subdomain SBP operators together into the global block diagonal matrix,
\begin{equation*}
\bar{\mathcal{D}}_{x_k}=  \sum_{m=1}^M (\tilde{I}_m^v)^\top \mathcal{D}_{x_k}^m\tilde{I}_m^v
= \begin{pmatrix}
\mathcal{D}_{x_k}^1 \\ & \ddots \\ && \mathcal{D}_{x_k}^M
\end{pmatrix}.
\end{equation*}
Note that $ \bar{\mathcal{D}}_{x_k}$ above defines the uncoupled (naive) $x_k-$derivative operator given by $ \bar{\mathcal{D}}_{x_k}=\mathcal{P}^{-1}\bar{\mathcal{Q}}_{x_k}$, where $\mathcal{P}$ and $\bar{\mathcal{Q}}_{x_k}$ possess the same block diagonal structure, i.e,
\begin{equation}\label{eq:SBP_blockdiag}
\mathcal{P} = \sum_{m=1}^M (\tilde{I}_m^v)^\top \mathcal{P}^m\tilde{I}_m^v
, \quad \bar{\mathcal{Q}}_{x_k} 
= \sum_{m=1}^M (\tilde{I}_m^v)^\top \mathcal{Q}_{x_k}^m\tilde{I}_m^v.
\end{equation}
Thus,  $\mathcal{P}$ in (\ref{eq:SBP_blockdiag}) defines a quadrature rule on the full set of volume nodes. 
With the notation from (\ref{eq:restr_ops}) and (\ref{eq:SBP_blockdiag}) above, we express the local SBP property (\ref{eq:SBP_split}) on the global mesh as follows,
\begin{equation}\label{eq:naive_SBP}
\begin{aligned}
\bar{ \mathcal{Q}}_{x_k}+\bar{\mathcal{Q}}_{x_k}^\top  = \ & 
(E^e)^\top P^eN_{x_k}^eE^e+ (E^i)^\top P^iN_{x_k}^iE^i,
\end{aligned}
\end{equation}
where we introduce the additional block matrices,
\begin{equation}\label{eq:SBP_blockdiag_boun}
\begin{aligned}
P ^e =  \ & \sum _{m=1}^M(\tilde{I}_m^e)^\top P_m^e  \tilde{I}_m^e ,  && N_{x_k} ^e= & \sum _{m=1}^M (\tilde{I}_m^e )^\top N_{x_k,m}^e  \tilde{I}_m^e ,  &&
E^e = \ & \sum _{m=1}^M E_m^e  \tilde{I}_m^v \\
P^i =  \ & \sum _{m=1}^M (\tilde{I}_m^v)^\top P_m^i \tilde{I}_m^v  && N_{x_k}^i = & \sum _{m=1}^M (\tilde{I}_m^i )^\top N_{x_k,m}^i  \tilde{I}_m^i ,  &&
E^i = \ & \sum _{m=1}^M E_m^i \tilde{I}_m^v,
\end{aligned}
\end{equation}
corresponding to the full set of external boundary and internal interface nodes, respectively.

For the approximation to satisfy the encapsulated SBP property (\ref{eq:SBP_def}) with respect to $P=P^e$, $N_{x_k}=N_{x_k}^e$ and $E=E^e$ in (\ref{eq:SBP_blockdiag_boun}), a penalty correction is required to cancel the interface contribution in (\ref{eq:naive_SBP}).
For this purpose, consider the following undivided difference matrix ansatz,
\begin{equation}\label{eq:SBP_multiblock}
\mathcal{Q}_{x_k} =\bar{\mathcal{Q}}_{x_k} - \frac{1}{2}(E^i)^\top P^i(N_{x_k}^i-N_{x_k}^*)E^i,
\end{equation}
where $\bar{\mathcal{Q}}_{x_k}$ is defined as in (\ref{eq:SBP_blockdiag}), and $N_{x_k}^*$ approximates $N_{x_k}^i$ via interface projection or interpolation operators.  Thus, (\ref{eq:SBP_multiblock}) introduces weak couplings between adjacent blocks or elements. From the SBP property on the global mesh (\ref{eq:naive_SBP}), the ansatz in (\ref{eq:SBP_multiblock}) yields
\begin{equation}\label{eq:SBP_multiblock_sufficient}
\mathcal{Q}_{x_k}+\mathcal{Q}_{x_k}^\top  = (E^e)^\top P^eN_{x_k}^eE^e+\frac{1}{2}(E^i)^\top \Big[P^iN_{x_k}^*+\big(P^iN_{x_k}^*\big)^\top \Big]E^i.
\end{equation}
Thus, a necessary and sufficient condition for (\ref{eq:SBP_def}) to hold with $P=P^e$, $N_{x_k}=N_{x_k}^e$ and $E=E^e$ is that the second term on the right hand side above is zero. Thus, $P^iN_{x_k}^*$ must be a skew-symmetric matrix.

For each $m\in\{1,\ldots, M\}$, let now $\mathcal{N}_m$ contain the index of all adjacent blocks or elements to $\Omega_m$ with at least partially overlapping faces. That is, we say that $n\in\mathcal{N}_m$ if $\int_{ \partial^i \Omega_n\cap \partial^i \Omega_m} \neq 0 $.  Next, we specify $N_{x_k}^*$ in (\ref{eq:SBP_multiblock}) to possess the off-diagonal block form given by
\begin{equation}\label{eq:SBP_multiblock_corr}
N_{x_k}^* =\sum_{m=1}^M\sum_{n\in\mathcal{N}_m}\frac{1}{2} (\tilde{I}_m^i )^\top \Big[ N_{x_k,m}^i  I_{nm}  + I_{nm}(-N_{x_k,n}^i ) \Big]\tilde{I}_n^i ,
\end{equation}
where $ I_{nm}$ denotes a discrete projection or interpolation operator from $\partial^i \Omega_n$ to $ \partial^i \Omega_m$. 
The pairwise inner product preserving (IPP) condition 
\begin{equation}\label{eq:IPP_def}
P_n^i  I_{mn}=(P_m^i I_{nm})^\top ,
\end{equation}
is also assumed to hold.
As the normal vectors point in opposite directions from the two opposing faces,  the summand in (\ref{eq:SBP_multiblock_corr}) should be understood as a mean value approximation rather than as a difference.

The following proposition generalizes the main stability result previously derived in \cite{Lundquist22_encaps} for the one-to-one coupling case.
\begin{proposition}
Consider the global multi-block/element $x_k-$derivative operator $\mathcal{D}_{x_k}=\mathcal{P}^{-1}\mathcal{Q}_{x_k}$ defined using (\ref{eq:SBP_blockdiag}), (\ref{eq:SBP_multiblock}) and (\ref{eq:SBP_multiblock_corr}). Provided that the interface operators satisfy the pairwise IPP condition (\ref{eq:IPP_def}), then the encapsulated SBP property (\ref{eq:SBP_def}) holds with respect to the full nodal mesh on $\Omega $.
\begin{proof}
From the definition of $N_{x_k}^*$ in (\ref{eq:SBP_multiblock_corr}), we have
\begin{equation*}
P^iN_{x_k}^* + (P^iN_{x_k}^*)^\top  = \frac{1}{2}\sum_{m=1}^M\sum_{n\in\mathcal{N}_m} R^{m,n} ,
\end{equation*}
where the summand $R^{m,n}$ is given by
 \begin{equation*}
 R^{m,n} = (\tilde{I}_m^i )^\top  P^{{m,i}}\big( N_{x_k,m}^i  I_{nm}  - I_{nm}N_{x_k,n}^i  \big)\tilde{I}_n^i  +(\tilde{I}_n^i )^\top \big(  I_{nm}^\top N_{x_k,m}^i   -N_{x_k,n}^i  I_{nm}^\top  \big)P_m^i \tilde{I}_m^i .
 \end{equation*}
Swapping the order of summation between $m$ and $n$ on the second term, we rewrite this sum as 
\begin{equation*}
\sum_{m=1}^M\sum_{n\in\mathcal{N}_m} R^{m,n} = \sum_{m=1}^M\sum_{n\in\mathcal{N}_m} \tilde{R}^{m,n} ,
\end{equation*}
where each term $\tilde{R}^{m,n}$ now comprises the single off-diagonal block term given by
 \begin{equation*}
 \begin{aligned}
\tilde{R}^{m,n} 
= \ &	(\tilde{I}_m^i )^\top  \Big[
P_m^i \big(N_{x_k,m}^i  I_{nm}  - I_{nm}N_{x_k,n}^i  \big)+ 
\big(  I_{mn}^\top N_{x_k,n}^i   -N_{x_k,m}^i  I_{mn}^\top  \big)P_n^i 
\Big]\tilde{I}_n^i .
\end{aligned}
 \end{equation*}
 Since $P_m^i $, $P_n^i $, $N_{x_k,m}^i $, $N_{x_k,n}^i $ are diagonal matrices,  they commute with respect to multiplication. So, we rearrange this expression to yield
 \begin{equation*}
\tilde{R}^{m,n} =(\tilde{I}_m^i )^\top  \Big[
N_{x_k,m}^i \big(P_m^i  I_{nm}  - I_{mn}^\top P_n^i  \big)+ 
\big(  I_{mn}^\top P_n^i   -P_m^i I_{nm} \big)N_{x_k,n}^i 
\Big]\tilde{I}_n^i .
\end{equation*}
Inserting the IPP condition (\ref{eq:IPP_def}) we now find that $\tilde{R}^{m,n}$ vanishes, i.e, $P^iN_{x_k}^*$ is a skew-symmetric matrix.  
In light of (\ref{eq:SBP_multiblock_sufficient}), this proves that $\mathcal{Q}_{x_k}$  satisfies (\ref{eq:SBP_def}).
\end{proof}
\end{proposition}
The practical construction of IPP  interface operators (\ref{eq:IPP_def}) for different types of quadrature rules is covered in the previous literature, including but not limited to  \cite{matt_carp10,Wang16_Tjunction,kozdon155,friedrich18, Lundquist18}.  Just as for SBP operators, tensor products and curvilinear transformations can be applied to extend IPP operators from $1$D to handle curved, non-conforming interfaces in multiple space dimensions.
For example, if curvilinear coordinates are used, then the IPP condition is preserved using the following rescalings,
\begin{equation}\label{eq:IPP_scaled}
I_{mn}=\sqrt{\underline{\boldsymbol{j}_n^i}}^{-1} \hat{I}_{mn}\sqrt{\underline{\boldsymbol{j}_m^i}} \quad I_{nm}=\sqrt{\underline{\boldsymbol{j}_m^i}}^{-1} \hat{I}_{nm}\sqrt{\underline{\boldsymbol{j}_n^i}},
\end{equation}
see \cite{Lundquist18} for more details.
Indeed, it is easily be verified that if the operators $\hat{I}_{mn}$ and $\hat{I}_{mn}$ are IPP with respect to $\hat{P}_m^i$ and $\hat{P}_n^i$, then $I_{mn}$ and $I_{mn}$ as defined above are IPP with respect to $P_n^i=\underline{\boldsymbol{j}_m^i}\hat{P}_m^i$ and $P_n^i=\underline{\boldsymbol{j}_n^i}\hat{P}_n^i$, see again (\ref{eq:curvi_boun_P}). Unfortunately, the exact consistency of $\hat{I}_{mn}$ and $\hat{I}_{mn}$ is not preserved by (\ref{eq:IPP_scaled}). Thus, if exact free-stream preservation is a desired property of the scheme, then a bit more care is needed, see e.g. \cite{Lundquist20_JCP} for complete details.

\subsection{Matrix-free implementation}
In a matrix-free implementation environment, the computation of matrix-vector products involving encapsulated multi-block/element SBP operators can be realized in different ways.
We start by splitting the expression for $\mathcal{Q}_{x_k}$ in (\ref{eq:SBP_multiblock}),  (\ref{eq:SBP_multiblock_corr}) into one block diagonal and one skew-symmetric off-diagonal part as follows,
\begin{equation}\label{eq:SBP_multiblock_alt}
\mathcal{Q}_{x_k} = \bar{\bar{\mathcal{Q}}}_{x_k}+ \frac{1}{2}(E^i)^\top P^i N_{x_k}^* E^i,
\end{equation}
where the block diagonal part is 
\begin{equation}\label{eq:Q_bar_bar}
\bar{\bar{\mathcal{Q}}}_{x_k}  =  \bar{\mathcal{Q}}_{x_k} - \frac{1}{2}(E^i)^\top P^i N_{x_k}^i E^i.
\end{equation}
Furthermore, the coupling term $N_{x_k}^*$ defined in (\ref{eq:SBP_multiblock_corr}) is compactly given by
\begin{equation}\label{eq:SBP_multiblock_corr_alt}
N_{x_k}^* =\frac{1}{2}(N_{x_k}^iI^i -I^iN_{x_k}^i),  \quad \text{where } I^i = \sum_{m=1}^M\sum_{n\in\mathcal{N}_m} (\tilde{I}_m^i )^\top I_{nm} \tilde{I}_n^i .
\end{equation}
Since the non-zero pattern of $(E^i)^TP^i N_{x_k}^i E^i$ in (\ref{eq:Q_bar_bar}) is either identical to or a subset of the non-zero pattern of $\bar{\mathcal{Q}}_{x_k} $, it makes practical sense to precompute and store the value of this block diagonal contribution $\bar{\bar{\mathcal{Q}}}_{x_k}$.  
In contrast, the off-diagonal contribution in (\ref{eq:SBP_multiblock_alt}) resembles an outer product in the sense that, in the generalized SBP case, the matrix-free operation $(E^i)^\top \big[P^i N_{x_k}^* (E^i\boldsymbol{U})\big]$ is cheaper to evaluate than the corresponding full matrix operation $\big[(E^i)^\top P^i N_{x_k}^* E^i\big]\boldsymbol{U}$. Therefore, although very useful for the purpose of analysis with the energy method,  it may not be the most computationally efficient approach to store and use the full matrix $\mathcal{Q}_{x_k}$.  Instead, the form (\ref{eq:SBP_multiblock_alt}) above provides a more suitable template for efficient matrix-free implementation.

We finally note that if curvilinear coordinates are employed in each subdomain according to the discussion in Section \ref{sec:curvi}, then storing the matrix $\bar{\bar{\mathcal{Q}}}_{x_k}$ in (\ref{eq:Q_bar_bar}) is suboptimal in terms of memory efficiency. Instead, we can make explicit use of the reference domain operators and metric coefficient vectors.  Starting from the definitions in (\ref{eq:curvi_disc_x}) through (\ref{eq:curvi_boun}), the block diagonal matrix $\bar{\bar{\mathcal{Q}}}_{x_k}$ in (\ref{eq:Q_bar_bar}) can then be rewritten as
\begin{equation}\label{eq:Q_bar_bar_curvi}
\begin{aligned}
\bar{\bar{\mathcal{Q}}}_{x} =  \ & \bar{\mathcal{Q}}_{x} -\frac{1}{2}(E^i)^\top \hat{P}^i  \big(N_\xi^i E^i\underline{\boldsymbol{Y}_\eta} -N_\eta^i E^i\underline{\boldsymbol{Y}_\xi} \big)  
\\
\bar{\bar{\mathcal{Q}}}_{y} =  \ & \bar{\mathcal{Q}}_{y} -\frac{1}{2}(E^i)^\top \hat{P}^i \big(N_\eta^i E^i\underline{\boldsymbol{X}_\xi} -N_\xi^i E^i\underline{\boldsymbol{X}_\eta} \big),
\end{aligned}
\end{equation}
where the block diagonal parts $\bar{\mathcal{Q}}_{x}$ and $\bar{\mathcal{Q}}_{y}$ satisfy the format from (\ref{eq:curvi_disc_x}), i.e,
\begin{equation*}
\begin{aligned}
\bar{\mathcal{Q}}_{x} =  \ & \tilde{\bar{\mathcal{Q}}}_x+\frac{1}{2}(E^e)^\top \hat{P}^e\Big[N_\xi^e\big(\underline{\boldsymbol{y}_\eta^e}E^e-E^e\underline{\boldsymbol{Y}_\eta}\big)-N_{\eta}^e \big(\underline{\boldsymbol{y}_\xi^e}E^e-E^e\underline{\boldsymbol{Y}_\xi}\big)\Big] 
\\
\bar{\mathcal{Q}}_{y} =  \ &\tilde{\bar{\mathcal{Q}}}_y+\frac{1}{2}(E^e)^\top \hat{P}^e\Big[N_\eta^e\big(\underline{\boldsymbol{x}_\xi^e}E^e-E^e\underline{\boldsymbol{X}_\xi}\big)-N_{\xi}^e \big(\underline{\boldsymbol{x}_\eta^e}E^e-E^e\underline{\boldsymbol{X}_\eta}\big)\Big] .
\end{aligned}
\end{equation*}
Finally,
 $\tilde{\bar{\mathcal{Q}}}_x$ and $ \tilde{\bar{\mathcal{Q}}}_y$ are given by 
\begin{equation*}
\begin{aligned}
\tilde{\bar{\mathcal{Q}}}_x= \ & \frac{1}{2} \big(\bar{\mathcal{Q}}_\xi\underline{\boldsymbol{Y}_\eta} +\underline{\boldsymbol{Y}_\eta}\bar{\mathcal{Q}}_\xi -\bar{\mathcal{Q}}_{\eta}\underline{\boldsymbol{Y}_\xi}-\underline{\boldsymbol{Y}_\xi}\bar{\mathcal{Q}}_{\eta}\big) \\
\tilde{\bar{\mathcal{Q}}}_y= \ & \frac{1}{2}\big(\bar{\mathcal{Q}}_{\eta}\underline{\boldsymbol{X}_\xi}+\underline{\boldsymbol{X}_\xi}\bar{\mathcal{Q}}_{\eta}- \bar{\mathcal{Q}}_\xi\underline{\boldsymbol{X}_\eta} - \underline{\boldsymbol{X}_\eta}\bar{\mathcal{Q}}_\xi \big),
 \end{aligned}
\end{equation*}
i.e, the same format as (\ref{eq:curvi_Q_tilde}).

\section{Accuracy}
\label{sec:accuracy}
As mentioned previously, a suboptimal accuracy bound related to classical SBP schemes for non-conforming block or element interfaces was proven in \cite{Lundquist18}.  In particular, this result applies to the auxiliary IPP interface operators (\ref{eq:IPP_def}) required for a stable implementation of the SBP operators.
For this type of mesh, fewer degrees of freedom are therefore required for the same order of approximation when employing generalized as opposed to classical SBP operators of a given accuracy order.  On the other hand, as we have seen, additional computational work is introduced in the form of outer products (\ref{eq:SBP_multiblock_alt}), the cost of which are proportional to the number of interface nodes as well as to the projection stencil widths. 
In practice, the overall computational efficiency also depends on
implementation details such as parallel computing models, the degree of parallelism as well as computer hardware. 
It is outside the scope of this work to provide detailed comparisons of efficiency between classical and generalized SBP operator implementations.
Instead, we present a discussion on the more objective metric given by operator accuracy.
In particular, we make an explicit comparison of accuracy between classical and generalized SBP approximations of collocated pseudo-spectral type as discussed Section \ref{sec:poly}.

For simplicity, consider two adjacent affine elements  $\Omega_m$ and $\Omega_n$ with partially overlapping faces. Without loss of generality, in the $2$D case we assume that this element interface is aligned in the $y-$direction.  For notational convenience we drop the superscript $i$ from all discrete interface quantities in the following discussion, and we thus let $\boldsymbol{y}_m$,  $\boldsymbol{y}_n$ contain the collocation quadrature nodes on the two opposing element faces, and $P_m $, $P_n $ the associated quadrature weights.
To define a set of pairwise of IPP interface operators (\ref{eq:IPP_def}), we employ $L_2$ projections \cite{Mavriplis89,kop96,kop02,friedrich18} in the following way,
\begin{equation}\label{eq:proj_dG}
 I_{nm} = P_m^{-1}\int_{\partial\Omega_n\cap \partial \Omega_m}\mathcal{L}_m \mathcal{L}_n^\top  dy, \quad   I_{mn} = P_n^{-1}\int_{\partial\Omega_m\cap \partial \Omega_n}\mathcal{L}_n \mathcal{L}_m^\top  dy, 
\end{equation}
where $\mathcal{L}_m=\mathcal{L}_m (y)$ and $\mathcal{L}_n= \mathcal{L}_n(y)$ contain the Lagrange polynomial functions on $\boldsymbol{y}_m$ and $\boldsymbol{y}_n$, respectively. From (\ref{eq:proj_dG}), we directly have
\begin{equation*}
P_m  I_{nm} = \int_{\partial\Omega_n\cap \partial \Omega_m} \mathcal{L}_m \mathcal{L}_n^\top dy=\Big(\int_{\partial\Omega_n\cap \partial \Omega_m}\mathcal{L}_n \mathcal{L}_m^\top dy\Big)^\top  = (P_n I_{mn})^\top ,
\end{equation*}
i.e, the IPP condition (\ref{eq:IPP_def}) is satisfied by construction.

As before, let $\mathcal{N}_m$ contain the index of all adjacent elements sharing the same partially overlapping face with $\Omega_m$. The overall degree of interpolation accuracy $p$ is given by the number of grid monomials which are projected exactly, i.e,
\begin{equation}\label{eq:dG_acc}
\sum_{n\in \mathcal{N}_m} I_{nm}\boldsymbol{y}_n^j = \boldsymbol{y}_m^j, \quad j = 1,\ldots ,  p,
\end{equation}
where the superscript $j$ indicates elementwise exponentiation.
Let $N_m+1$ and $N_n+1$ be the number of quadrature nodes, i.e. the length of $\boldsymbol{y}_m$ and $\boldsymbol{y}_n$, respectively.  By definition, the interpolating Lagrange polynomials satisfy
\begin{equation*}
\mathcal{L}_n^\top \boldsymbol{y}_n^j =y^j  , \quad j \leq N_n.
\end{equation*}
We thus have
\begin{equation*}
P_m I_{nm}\boldsymbol{y}_n^j = \int_{\partial\Omega_n\cap \partial \Omega_m} \mathcal{L}_m (\mathcal{L}_n^\top \boldsymbol{y}_n^j)dy = \int_{\partial \Omega_n\cap \partial  \Omega_m } \mathcal{L}_my^j  dy, \quad j \leq N_n.
\end{equation*}
Summing over all adjacent faces, this yields
\begin{equation*}
P_m \sum_{n\in \mathcal{N}_m} I_{nm}\boldsymbol{y}_n^j =\sum_{n\in \mathcal{N}_m} \int_{\partial\Omega_n\cap \partial \Omega_m}  \mathcal{L}_m y^j dy = \int_{\partial  \Omega_m}  \mathcal{L}_m y^jdy, \quad j \leq \underset{n\in\mathcal{N}_m}{\mathrm{min}}(N_n).
\end{equation*}
Now, let $q_m$ denote the quadrature degree of $P_m $. Since the vector $\mathcal{L}_m y^j$ contains polynomials of order $N_m  + j$, we can then further write
\begin{equation*}
\int_{\partial  \Omega_m}  \mathcal{L}_m y^j dy= \int_{P_m }  \mathcal{L}_m y^jdy, \quad j \leq q_m - N_m,
\end{equation*}
where $\int_{P_m }$ indicates numerical integration with the quadrature weights associated with the $N_m+1$ quadrature nodes. Finally, we have
\begin{equation*}
\int_{P_m }  \mathcal{L}_m y^j= \int_{P_m }  \mathcal{L}_m(\mathcal{L}_m^\top  \boldsymbol{y}_m^j )dy = \Big(\int_{P_m }  \mathcal{L}_m\mathcal{L}_m^\top   dy\Big)\boldsymbol{y}_m^j= P_m \boldsymbol{y}_m^j, \quad j\leq N_m,
\end{equation*}
where the last equality follows from the Kronecker delta property of Lagrange interpolating polynomials.

Combining the results above, we have shown that (\ref{eq:dG_acc}) is satisfied for $p= \mathit{min}(N_m, \ N_n,\ q_m-N_m)$.
In other words, the degree of $L_2$ projection accuracy is constrained by $q_m-N_m=N_m-1$ in the Gauss-Lobatto quadrature case, which is subobtimal compared to the polynomial degree of $N_m$. However, no such order reduction occurs if Gauss-Legendre quadrature nodes are used, as in this case we have $p=q_m-N_m=N_m+1$.  Recall that these two types of quadrature rules are associated with classical and generalized pseudo-spectral SBP operators, respectively. 
\begin{remark}
The result of a higher quadrature degree allowing for more accurate IPP interpolation is in line with a more general bound on IPP operator accuracy derived in \cite{Lundquist20_JCP}.
\end{remark}

\section{Numerical results}
In order to investigate the practical implications of a suboptimal $L_2$ projection accuracy $p$ in (\ref{eq:dG_acc}), in this section we compare the overall accuracy resulting from the encapsulated multi-domain SBP operator formulation in (\ref{eq:SBP_multiblock}), (\ref{eq:SBP_multiblock_corr}) using classical as well as generalized SBP operators of pseudo-spectral type.

\subsection{Operator accuracy}
\begin{figure}[htb!]
\begin{center}
    \includegraphics[width=.5\textwidth]{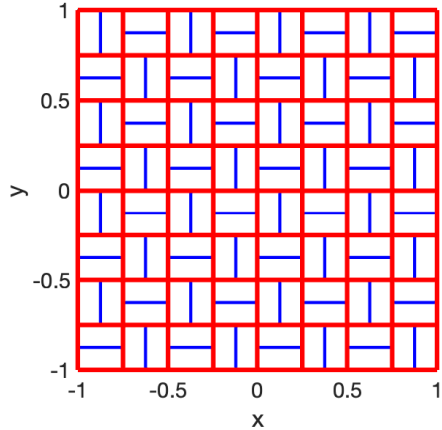}
    \end{center}    
   \caption{Non-conforming multi-element mesh structure.  Each rectangle represents a single $2$D pseudo-spectral element.}
    \label{fig:elem_mesh}
\end{figure}

We start by comparing the approximation errors resulting from applying the global SBP operators to a known, smooth function. For this test, we consider the multi-element mesh structure shown in Figure  \ref{fig:elem_mesh}. Each small rectangle is discretized using pseudo-spectral SBP operators (\ref{eq:pseudo-op}) with a given polynomial degree $N$, and the two-to-one element couplings are facilitated using IPP interface operators of $L_2$ projection type (\ref{eq:proj_dG}).  Since the mesh is non-conforming at every element interface throughout the computational domain, the resulting global operators are especially sensitive to large $L_2$ projection errors.

\begin{figure}[htb!]
  \begin{subfigure}{0.5\textwidth}
    \includegraphics[width=\linewidth]{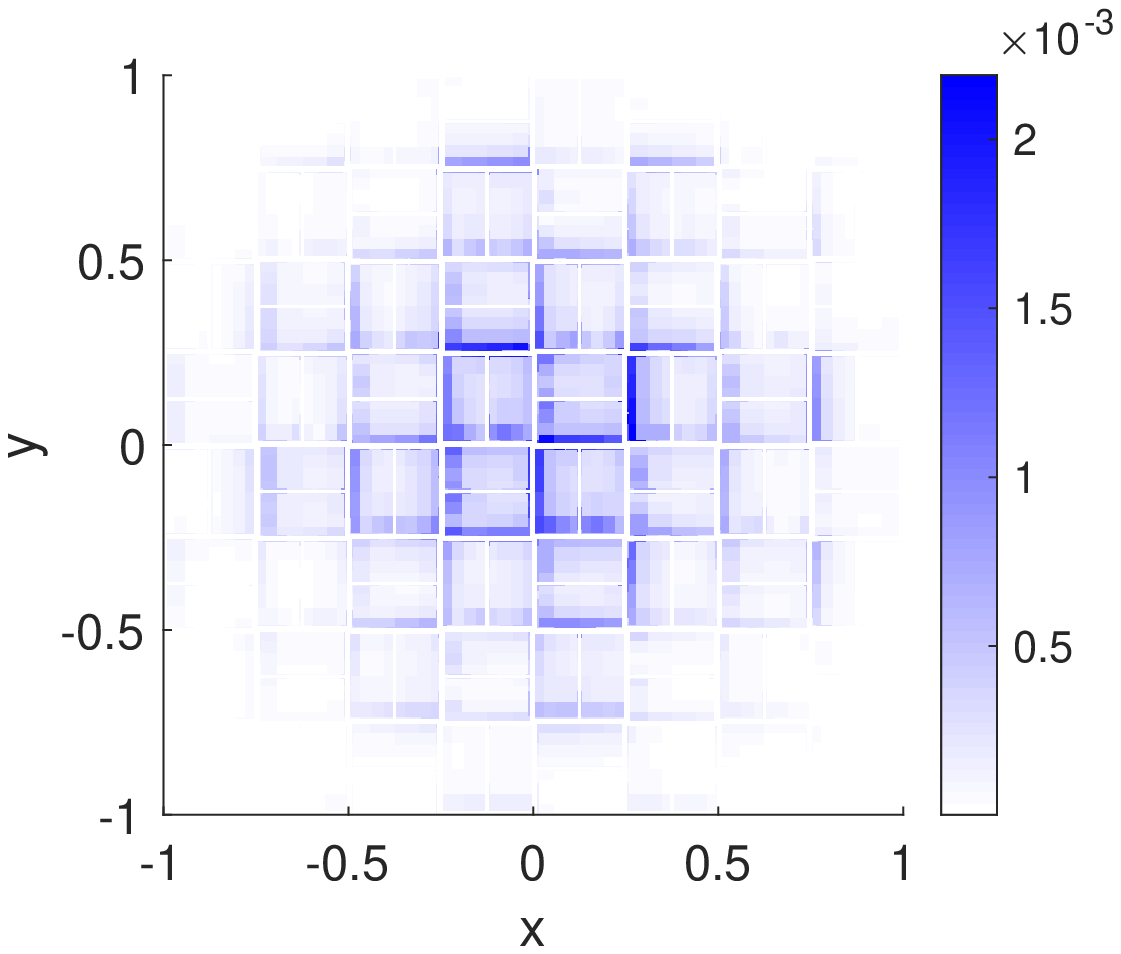}
   \caption{Gauss-Legendre} \label{fig:LG4_err}
  \end{subfigure}%
  \hfill   
  \begin{subfigure}{0.5\textwidth}
    \includegraphics[width=\linewidth]{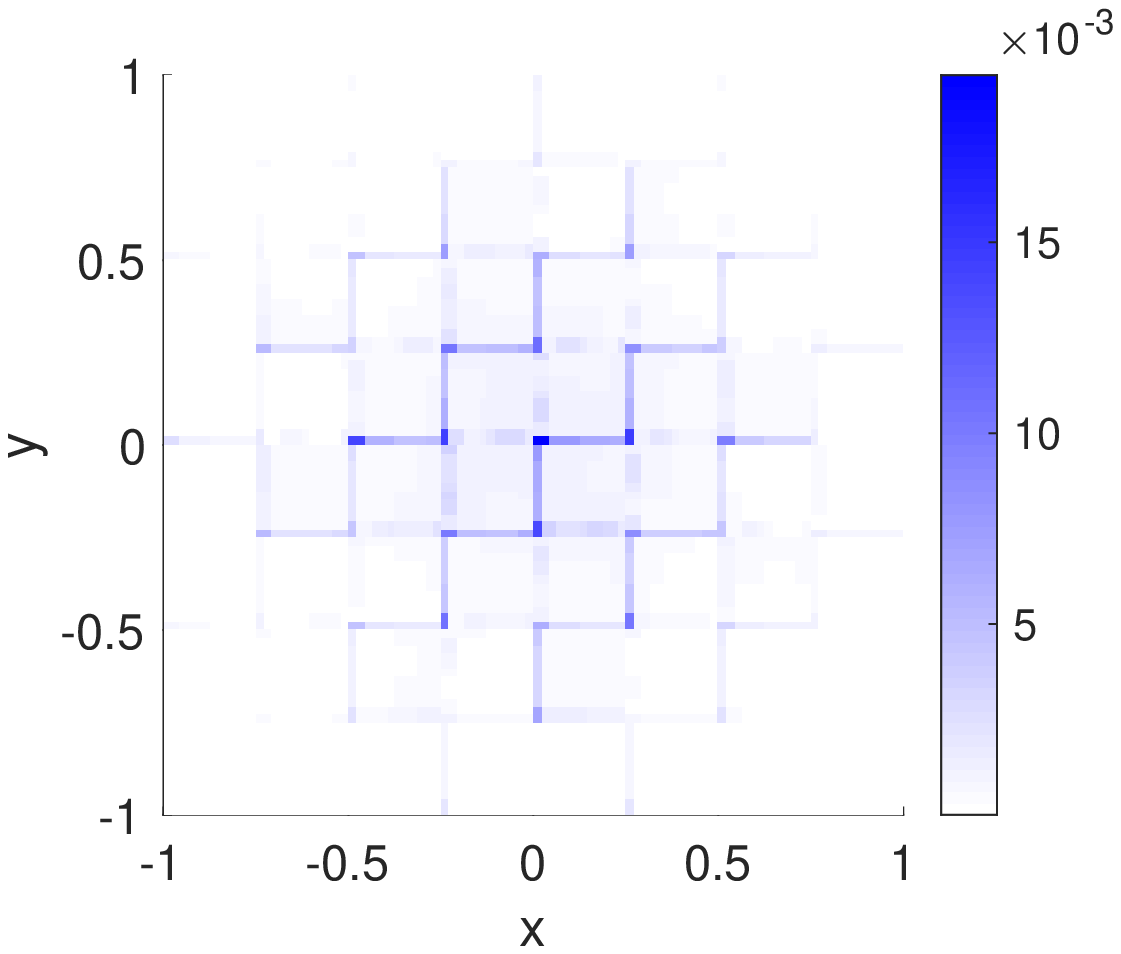}
    \caption{Gauss-Lobatto}  \label{fig:LGL4_err}
  \end{subfigure}%
\caption{Absolute error levels using polynomial degree $N=4$. Note the different color bar scales.} \label{fig:LG_errors}
\end{figure}

\begin{table*}[htb!]\centering
\begin{tabular}{@{}ccccccccccc@{}}
& \multicolumn{2}{c}{Gauss-Legendre} & \phantom{abc}& \multicolumn{2}{c}{Gauss-Lobatto}\\
 \cmidrule{2-3} \cmidrule{5-6} \cmidrule{8-9}
$n$ & $|\boldsymbol{e}|_{\mathrm{ max}}$ & $\text{rate}$ && $|\boldsymbol{e}|_{\mathrm{ max}}$ & $\text{rate}$
\\ \midrule
     1  & 2.0098e-02  &    -      &&        1.1158e-01      &    -    &&      \\
    2 & 3.0909e-03  &  2.70  && 3.1727e-02  &    1.81    &&    \\
    4 & 4.0608e-04&  2.93  && 8.1558e-03  & 1.96  &&  \\
   8 & 5.1390e-05  &  2.98  && 2.0489e-03 & 1.99  &&  \\
        16 & 6.4436e-06  &  3.00  && 5.1231e-04   & 2.00  &&  \\
 \bottomrule
\end{tabular}
\caption{Error convergence using polynomial degree $N=3$.}
\label{table:conv3} 
\end{table*}

\begin{table*}[htb!]\centering
\begin{tabular}{@{}ccccccccccc@{}}
& \multicolumn{2}{c}{Gauss-Legendre} & \phantom{abc}& \multicolumn{2}{c}{Gauss-Lobatto}\\
 \cmidrule{2-3} \cmidrule{5-6} \cmidrule{8-9}
$n$ & $|\boldsymbol{e}|_{\mathrm{ max}}$ & $\text{rate}$ && $|\boldsymbol{e}|_{\mathrm{ max}}$ & $\text{rate}$
\\ \midrule
     1  & 2.1900e-03  &    -      &&        1.9378e-02      &    -    &&      \\
    2 & 1.5270e-04  &  3.84  && 2.9354e-03  &    2.72    &&    \\
    4 & 9.5342e-06 &  4.00  && 3.8473e-04  & 2.93  &&  \\
   8 & 6.0808e-07  &  3.97  && 4.8663e-05 & 2.98  &&  \\
        16 & 3.8038e-08  &  4.00  && 6.1008e-06   & 3.00  &&  \\
 \bottomrule
\end{tabular}
\caption{Error convergence using polynomial degree $N=4$.}
\label{table:conv4} 
\end{table*}

\begin{table*}[htb!]\centering
\begin{tabular}{@{}ccccccccccc@{}}
& \multicolumn{2}{c}{Gauss-Legendre} & \phantom{abc}& \multicolumn{2}{c}{Gauss-Lobatto}\\
 \cmidrule{2-3} \cmidrule{5-6} \cmidrule{8-9}
$n$ & $|\boldsymbol{e}|_{\mathrm{ max}}$ & $\text{rate}$ && $|\boldsymbol{e}|_{\mathrm{ max}}$ & $\text{rate}$
\\ \midrule
     1  & 2.3645e-04  &    -      &&        2.6002e-03      &    -    &&      \\
    2 & 1.0380e-05  &  4.51  && 1.8397e-04  &    3.82    &&    \\
    4 & 3.5035e-07 &  4.89  && 1.1541e-05  & 3.99  &&  \\
   8 & 1.1155e-08  &  4.97  && 7.2441e-07 & 3.99  &&  \\
        16 & 3.4496e-10  &  5.00  && 4.5266e-08  & 4.00  &&  \\
 \bottomrule
\end{tabular}
\caption{Error convergence using polynomial degree $N=5$.}
\label{table:conv5} 
\end{table*}

We employ the exact smooth function given by $u=e^{-((3x)^2+(3y)^2)/2}$. Let $\boldsymbol{U}_x$ and $\boldsymbol{U}_y$ denote the exact partial derivatives in the volume nodes, and consider the pointwise error measure given by $|\boldsymbol{e}|=\sqrt{|\mathcal{D}_x\boldsymbol{U}-\boldsymbol{U}_x|^2+|\mathcal{D}_y\boldsymbol{U}-\boldsymbol{U}_y|^2}$. The numerical results for $N=4$ are shown in Figure \ref{fig:LG_errors} using both classical (Gauss-Lobatto) and generalized (Gauss-Legendre) SBP operators. Note that the number of mesh nodes is the same in both cases. As can be seen, the pointwise error levels are especially elevated at element interfaces in the classical SBP case. This result is explained by the suboptimal $L_2$ projection accuracy order of $p=N-1$ associated with the Gauss-Lobatto quadrature nodes.

Convergence studies for three polynomial degree values are given in Table \ref{table:conv3}, Table \ref{table:conv4} and Table \ref{table:conv5}.  The coarsest refinement level ($n=1$) here corresponds to the multi-element mesh shown in Figure \ref{fig:elem_mesh}. The convergence rates obtained are consistently one order less than the IPP interpolation accuracy order given by $p=N$ and $p=N-1$, respectively. 
This behaviour is explained by the combined scaling of $\mathcal{D}_{x_i}=\mathcal{P}^{-1}\mathcal{Q}_{x_i}$ with $\mathcal{P}^{-1}$ and $P^i$,  see (\ref{eq:SBP_multiblock}), where the volume quadrature weights in $\mathcal{P}$ are smaller than the interface quadrature weights in $P^{i}$ by one order of magnitude.
The results are thus in line with theoretical expectations, confirming that generalized SBP operators are associated with smaller overall approximation errors by one order compared to classical SBP operators for a given mesh size.

\subsection{Scalar advection problem}
\begin{figure}[htb!]
\begin{center}
    \includegraphics[width=.5\textwidth]{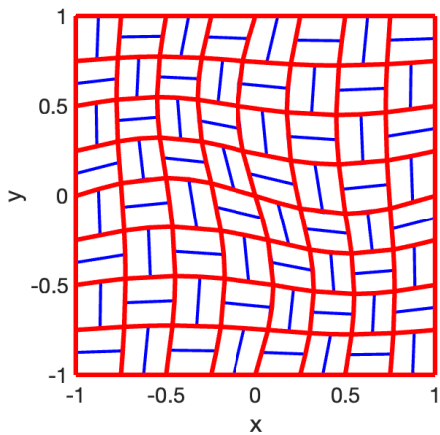}
    \end{center}    
   \caption{Curvilinear non-conforming mesh structure.}
    \label{fig:curvi_mesh}
\end{figure}

Next, we apply encapsulated SBP operators in order to semi-discretize the linear advection model problem
\begin{equation*}
u_t + u_x + u_y = 0, \quad 0 \leq t \leq 0.5.
\end{equation*}
For this test, we consider the curvilinear multi-element mesh structure illustrated in Figure \ref{fig:curvi_mesh}. Compared to the Cartesian mesh shown in Figure \ref{fig:elem_mesh}, we have here applied the following coordinate transformation,
\begin{equation*}
\begin{aligned}
x(\xi,  \eta) = \ & x + 0.1\sin \pi (\eta + 1) \\
y(\xi,  \eta) = \ & y + 0.1\sin \pi (\xi + 1),
\end{aligned}
\end{equation*}
and the explicit formulae (\ref{eq:SBP_multiblock_alt}),  (\ref{eq:SBP_multiblock_corr_alt}), (\ref{eq:Q_bar_bar_curvi}) are used in order to construct a pair of encapsulated SBP operators (\ref{eq:SBP_def}) on the full mesh.

\begin{figure}[htb!]
  \begin{subfigure}{0.5\textwidth}
    \includegraphics[width=\linewidth]{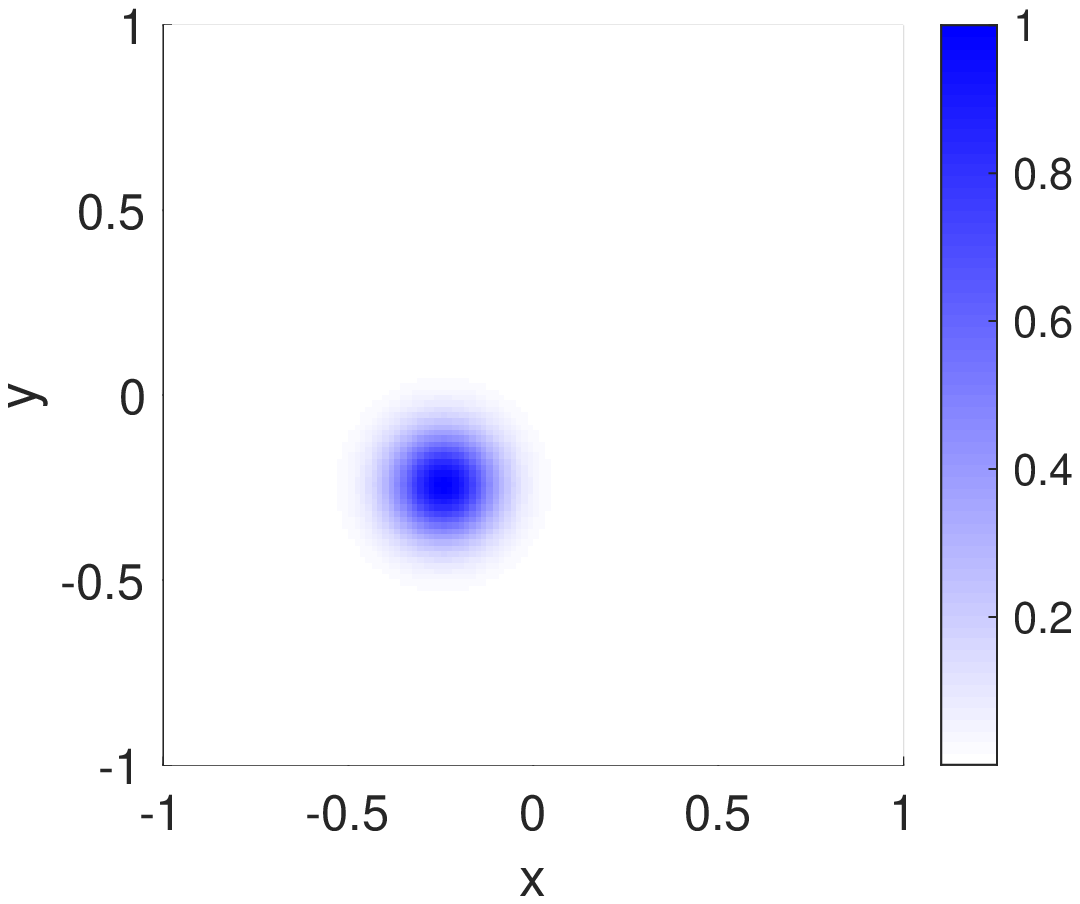}
   \caption{$t=0$} \label{fig:exact_t0}
  \end{subfigure}%
  \hfill   
  \begin{subfigure}{0.5\textwidth}
    \includegraphics[width=\linewidth]{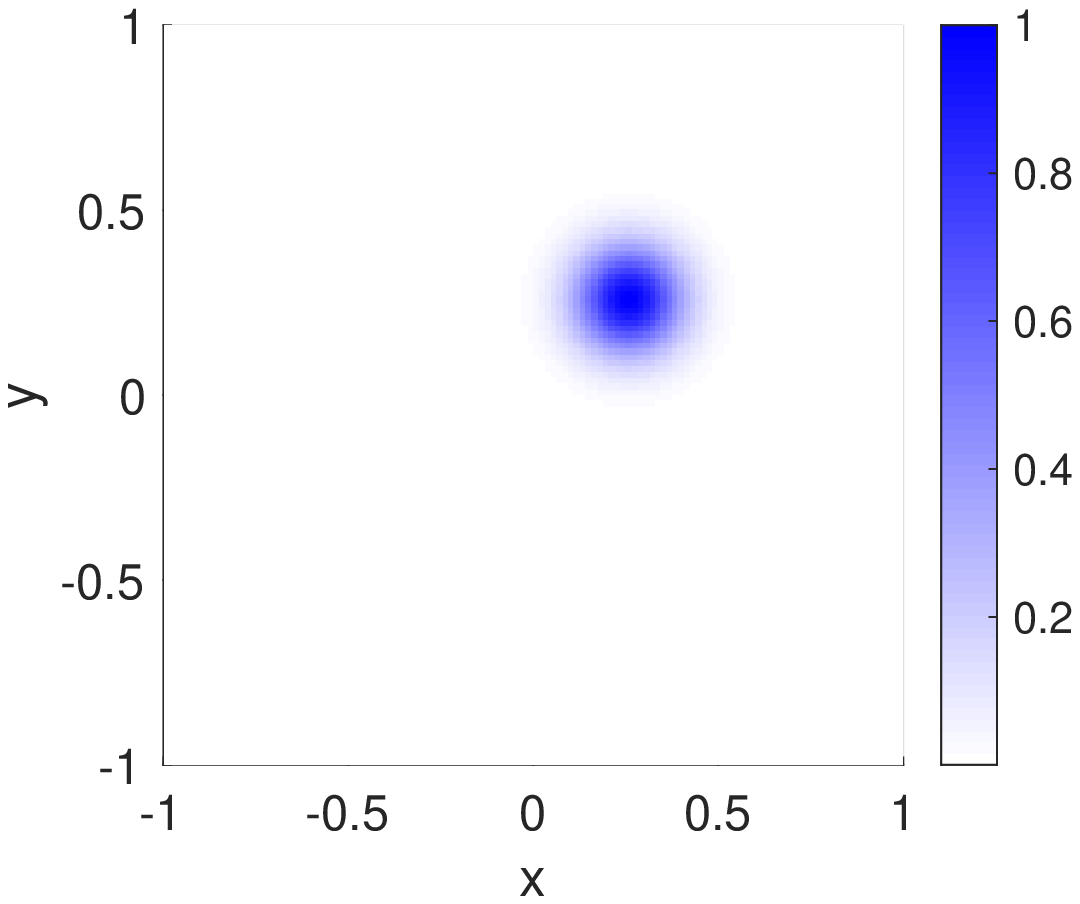}
    \caption{$t=0.5$}  \label{fig:exact_t1}
  \end{subfigure}%
\caption{Exact solution at the initial and final time.} \label{fig:exact_curvi}
\end{figure}
\begin{figure}[htb!]
  \begin{subfigure}{0.5\textwidth}
    \includegraphics[width=\linewidth]{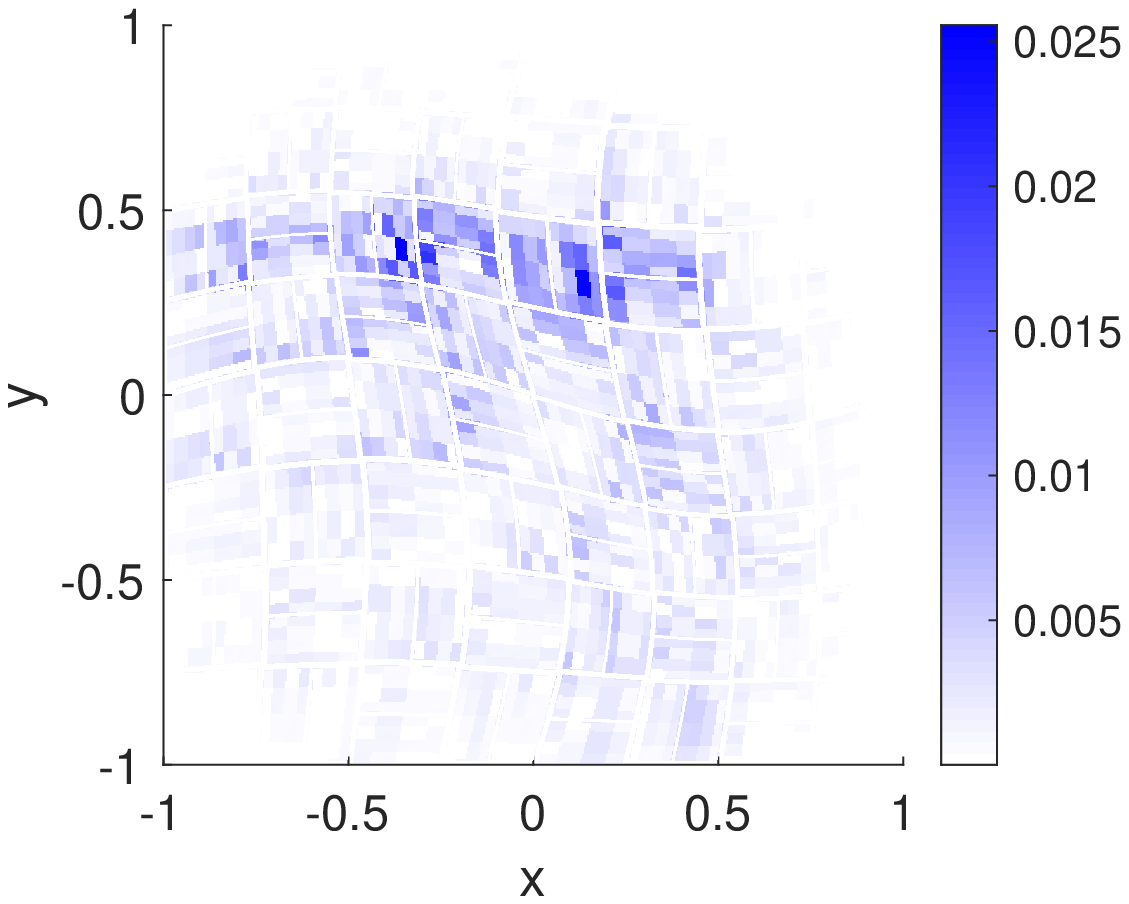}
   \caption{Gauss-Legendre} \label{fig:LG4_err_curvi}
  \end{subfigure}%
  \hfill   
  \begin{subfigure}{0.5\textwidth}
    \includegraphics[width=\linewidth]{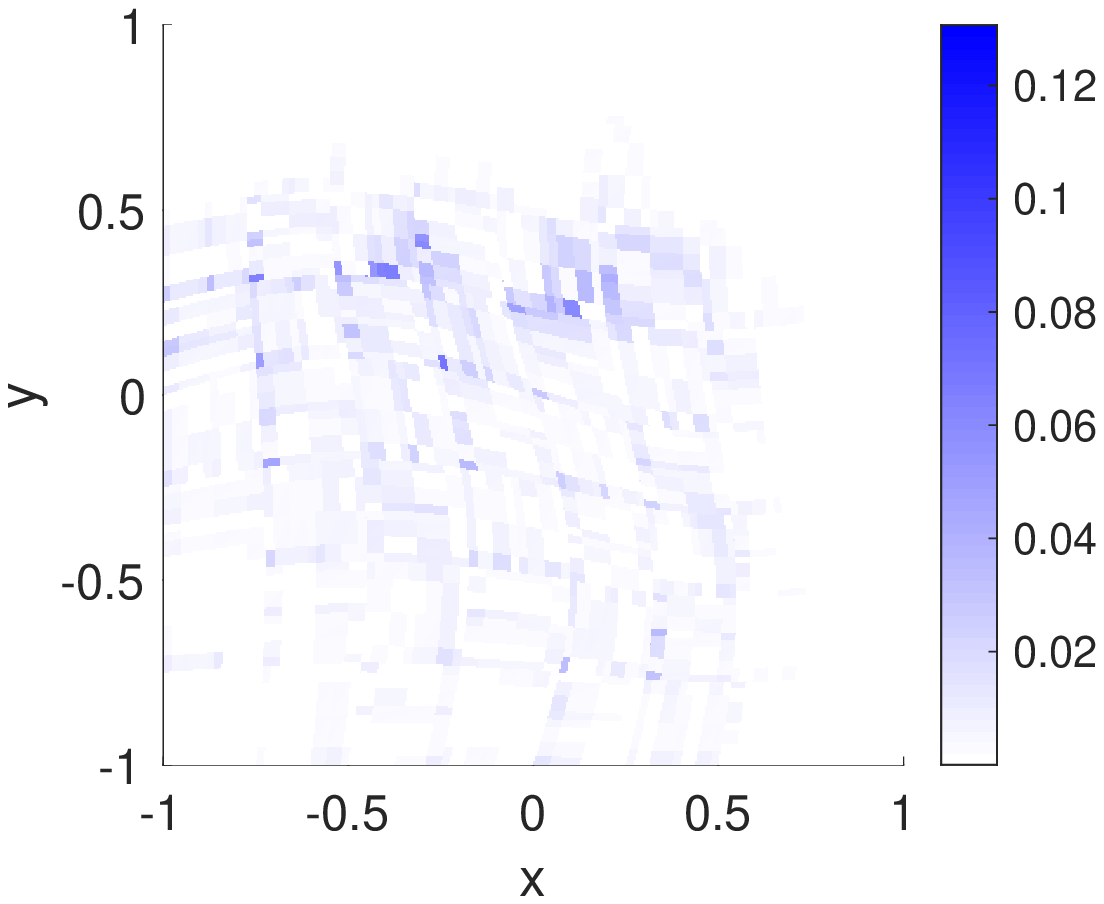}
    \caption{Gauss-Lobatto}  \label{fig:LGL4_err_curvi}
  \end{subfigure}%
\caption{Solution error at $t=0.5$ on cuvilinear mesh.  Note the different color bar scales.} \label{fig:LG_errors_curvi}
\end{figure}
\begin{table*}[htb!]\centering
\begin{tabular}{@{}ccccccccccc@{}}
& \multicolumn{2}{c}{Gauss-Legendre} & \phantom{abc}& \multicolumn{2}{c}{Gauss-Lobatto}\\
 \cmidrule{2-3} \cmidrule{5-6} \cmidrule{8-9}
$n$ & $|\boldsymbol{e}|_{\mathrm{ max}}$ & $\text{rate}$ && $|\boldsymbol{e}|_{\mathrm{ max}}$ & $\text{rate}$
\\ \midrule
     1  & 5.6588e-02  &    -      &&        4.0673e-01      &    -    &&      \\
    2 & 1.0215e-02  &  2.47  && 5.8035e-02  &    2.81    &&    \\
    4 & 6.6422e-04&  3.94  && 1.2944e-02  & 2.16  &&  \\
   8& 1.0689e-04  &  2.63  && 2.9054e-03 & 2.16  &&  \\
        16 & 2.8676e-05  &  1.90  && 4.3580e-04   & 2.74  &&  \\
         av. &   &  2.74  &&    & 2.47  &&  \\
 \bottomrule
\end{tabular}
\caption{SBP approximation convergence using polynomial degree $N=3$.}
\label{table:conv3_curvi} 
\end{table*}

\begin{table*}[htb!]\centering
\begin{tabular}{@{}ccccccccccc@{}}
& \multicolumn{2}{c}{Gauss-Legendre} & \phantom{abc}& \multicolumn{2}{c}{Gauss-Lobatto}\\
 \cmidrule{2-3} \cmidrule{5-6} \cmidrule{8-9}
$n$ & $|\boldsymbol{e}|_{\mathrm{ max}}$ & $\text{rate}$ && $|\boldsymbol{e}|_{\mathrm{ max}}$ & $\text{rate}$
\\ \midrule
     1  & 2.5565e-02  &    -      &&        1.3051e-01      &    -    &&      \\
    2 & 1.4661e-03  &  4.12  && 9.6445e-03  &    3.75    &&    \\
    4 & 2.1635e-04 &  2.76  && 9.4095e-04  & 3.36  &&  \\
  8 & 6.1528e-06  &  5.14  && 9.6140e-05 & 3.29  &&  \\
        16 & 1.9403e-07  &  4.99  && 7.7472e-06   & 3.63  &&  \\
                 av.  &   &  4.25  &&    & 3.51  &&  \\
 \bottomrule
\end{tabular}
\caption{SBP approximation convergence using polynomial degree $N=4$.}
\label{table:conv4_curvi} 
\end{table*}

\begin{table*}[htb!]\centering
\begin{tabular}{@{}ccccccccccc@{}}
& \multicolumn{2}{c}{Gauss-Legendre} & \phantom{abc}& \multicolumn{2}{c}{Gauss-Lobatto}\\
 \cmidrule{2-3} \cmidrule{5-6} \cmidrule{8-9}
$n$ & $|\boldsymbol{e}|_{\mathrm{ max}}$ & $\text{rate}$ && $|\boldsymbol{e}|_{\mathrm{ max}}$ & $\text{rate}$
\\ \midrule
     1  & 5.9775e-03  &    -      &&        3.9743e-02      &    -    &&      \\
    2 & 2.7492e-04 &  4.44  && 1.4142e-03  &   4.82    &&    \\
    4 & 4.1758e-06 &  6.04  && 6.6759e-05  & 4.40  &&  \\
   8 & 1.5996e-07  &  4.71  && 2.8231e-06 & 4.56  &&  \\
       16 & 5.1622e-09  &  4.95  && 2.8136e-07  & 3.33  &&  \\
                         av.  &   &  5.04  &&    & 4.27  &&  \\
 \bottomrule
\end{tabular}
\caption{SBP approximation convergence using polynomial degree $N=5$.}
\label{table:conv5_curvi} 
\end{table*}

We employ the exact solution given by
\begin{equation*}
u = e^{-0.005\big[(x+0.25-t)^2-(y+0.25-t)^2\big]}.
\end{equation*}
The solution at the initial and final time is also shown in Figure \ref{fig:exact_curvi}. 
The semi-discrete system is marched forward in time using the classical fourth order explicit Runge-Kutta method. In order to compare the accuracy resulting from different spatial discretizations, the time step sizes are chosen sufficiently small to make temporal discretization errors negligible.  The resulting numerical errors are plotted in Figure \ref{fig:LG_errors_curvi} for the case $N=4$,  and convergence studies for $N=3$, $N=4$ and $N=5$ are given in Table \ref{table:conv3_curvi}, Table \ref{table:conv4_curvi} and Table \ref{table:conv5_curvi}, respectively. Although the convergence results are no longer as clear-cut, they still confirm that generalized SBP operators lead to more accurate solutions than classical ones for a given polynomial degree. On average,  the observed convergence appears to approximately follow the rates $N$ and $N-1/2$, respectively.

\section{Conclusions}
We have successfully extended both the construction and stable implementation of encapsulated summation-by-parts operators for curvilinear and multi-block/element meshes from a classical to a generalized summation-by-parts operator definition.  By simplifying both the analysis and implementation, this development promotes the successful usage of generalized summation-by-parts operators in the context of non-linear problems, curvilinear coordinates as well as local mesh refinement. In particular, the favourable accuracy properties of generalized compared to classical summation-by-parts schemes in the presence of non-conforming nodal interfaces can thus be harnessed in an optimal way. 
By the proper introduction of correction penalty terms at physical boundaries,  we have shown that numerical stability can be ensured for general nonlinear and variable coefficient problems with the energy method.
The superior accuracy associated with generalized as opposed to classical SBP operators of pseudo-spectral element type was also confirmed both theoretically and experimentally. The overall efficiency of generalized operators (which require more work) compared to classical operators remains unknown however, and requires additional future investigations.
\section*{Acknowledgments}
Tomas Lundquist and Andrew Winters were supported by Vetenskapsr{\aa}det, Sweden [award no.~2020-03642 VR]. Jan Nordstr\"{o}m was supported by Vetenskapsr{\aa}det, Sweden [award no.~2021-05484 VR].

%

\bibliographystyle{model1-num-names}
\bibliography{References}







\end{document}